\documentclass{siamart}
\usepackage{bm}
\usepackage{epsfig,multirow,bbm}
\usepackage{amsmath,enumitem}
\usepackage{amssymb}
\usepackage{url,changebar,mathtools}
\usepackage{subfigure}
\usepackage{cite}   
\usepackage{wrapfig}
\usepackage{amsfonts}
\usepackage{graphicx}
\usepackage{epstopdf}
\ifpdf
  \DeclareGraphicsExtensions{.eps,.pdf,.png,.jpg}
\else
  \DeclareGraphicsExtensions{.eps}
\fi

\usepackage{algorithm}
\usepackage[noend]{algpseudocode}
\setlist[enumerate]{leftmargin=.5in}
\setlist[itemize]{leftmargin=.5in}

\newsiamremark{remark}{Remark}
\newsiamremark{hypothesis}{Hypothesis}
\crefname{hypothesis}{Hypothesis}{Hypotheses}
\newsiamthm{claim}{Claim}

\headers{Extreme excursion probabilities for dynamical systems}{V. Rao and M. Anitescu}

\title{{Efficient computation of extreme excursion probabilities for dynamical systems}\thanks{This material was based upon work
	supported by the U.S. Department of Energy, Office of Science,
	Office of Advanced Scientific Computing Research (ASCR) under
	Contract DE-AC02-06CH11347. We acknowledge partial NSF funding
	through awards FP061151-01-PR and CNS-1545046 to MA.}}

\author{Vishwas Rao\thanks{Argonne National Laboratory, Lemont, IL 
  (\email{vhebbur@anl.gov}). Corresponding author.}
  \and Mihai Anitescu\thanks{Argonne National Laboratory, Lemont, IL and the University of Chicago (\email{anitescu@mcs.anl.gov}).} }

\usepackage{amsopn}

\newcommand{\x}{\mathbf{x}}
\newcommand{\y}{\mathbf{y}}

\begin{document}

\maketitle

\begin{abstract}
We develop a novel computational method for evaluating the extreme excursion probabilities arising for random initialization of nonlinear dynamical systems. The method uses a Markov chain Monte Carlo or a Laplace approximation approach to construct a biasing distribution that in turn is used in an importance sampling procedure to estimate the extreme excursion probabilities. The prior and likelihood of the biasing distribution are obtained by using Rice's formula from excursion probability theory. We use Gaussian mixture biasing distributions and approximate the non-Gaussian initial excitation by the method of moments to circumvent the linearity and Gaussianity assumptions needed by excursion probability theory. We demonstrate the effectiveness of this computational framework for  nonlinear dynamical systems of up to 100 dimensions.   
\end{abstract}
      
      % REQUIRED
      
      % REQUIRED
\section{Motivation}\label{sec:intro} 
Computing the probability of extreme events is of central importance for dynamical systems that arise in natural phenomena such as climate, weather, oceanography {\color{blue} \cite{Easterling_2000, Easterling_2000A}}, and engineering systems such as structures \cite{Cornell_1968, Vrouwenvelder_2000} or power grids \cite{Lesieutre_2008}. Examples of consequential extreme events are rogue waves in the ocean \cite{Dysthe_2008}, hurricanes, tornadoes \cite{Ross_2003}, and power outages \cite{Atputharajah_2009}. In this work we are motivated by the increased concern of transient security in the presence of uncertain inertia, as identified by the North American Electric Reliability Corporation  in its most recent long-term reliability assessment \cite{nerc2017ltra}. The mathematical formulation is  of a dynamical system with parametric uncertainty, which is equivalent to initial condition uncertainty by adding the equations $\dot{p}=0$ to the ordinary differential equation, where $p$ are the parameters. The aim of the calculation in that case is to compute an extreme excursion probability: the odds that the transient due to a sudden malfunction exceeds prescribed safety limits. Since the reliability goal is that an average customer should experience only minutes of electricity interruption per year \cite{fischer2016german}, the target safe limit exceedance probabilities may be in the range of $10^{-4}$--$10^{-5}$. 

 Quantifying extreme excursion probabilities is of great importance  because of their  socioeconomic impact. Outcomes of interest reside in the tails of the probability distribution of the associated event space because of their low likelihood. To resolve the tails of these events, one has to evaluate multivariable integrals over complex domains. Because of the tiny mass and complex shape of the relevant likelihood level sets, standard quadrature, cubature, or sparse grid  methods cannot be  applied directly to evaluate these integrals. The most commonly used method is Monte Carlo simulation (MCS), which requires repeated samples of the underlying outcome. For such small probabilities, however, MCS exhibits a large variance relative to the probability to be computed, and thus it needs a large number of samples to produce results of acceptable accuracy. For example, estimating the odds of an extreme event, whose probability ends up being $10^{-3}$ for an underlying process that requires 10 minutes per numerical simulation, requires two years of serial computation for producing an estimate with a standard deviation of less than 10\% of the target value via MCS. Hence, alternative methods must be developed that are computationally efficient. 

In the rest of this section, we review the literature, provide an overview of our approach, and discuss its limitations and possible extensions. In \S \ref{sec:problem} we describe the rare-event problem in detail and revisit MCS and importance sampling (IS) methods. In \S \ref{sec:mcmc-based-sampling} we formulate the problem of estimating rare-event probability as a sequence of Bayesian inverse problems, in \S \ref{sec:BayesianInverseProblem} we discuss two well-known approaches to solve the Bayesian inverse problems, and in \S \ref{sec:ibd} we use the solutions of these Bayesian inverse problems to construct an importance biasing distribution (IBD). In \S \ref{sec:num_exp} we demonstrate the algorithm on two nonlinear dynamical systems of different sizes.  In \S \ref{sec:conc} we give concluding remarks .
\subsection{Literature review}Most  methods to compute the probabilities of rare events are a combination of MCS and IS methods. The key difference between such approaches lies in the proposal distribution for importance sampling. In what follows, we briefly discuss existing methods and the key ideas underpinning them.
\subsubsection{Monte Carlo simulation}
The MCS approach is one of the most robust methods for simulating rare events and estimating their probabilities. It was originally developed for solving problems in mathematical physics \cite{Metropolis_1949}. Since then,  the method has been used in a wide variety of applications, and it currently lies at the heart of all random sampling-based techniques \cite{Liu_2008B, Robert_2005B}. The main strength of MCS is that its rate of convergence does not depend on the  likelihood level set or its dimension. When evaluating excursion probabilities, the method primarily counts how many of the random samples exceed the given excursion level. Thus, in order to estimate a probability $p$, MCS needs a number of samples exceeding $\frac{1}{p}$, which for small probabilities makes its direct application impractical. 
%GAIL - that seems like a disadvantage; should you use a qualifying word like "When evaluating..., however, ...
\subsubsection{Importance sampling} IS methods belong to the class of variance reduction techniques that aim to estimate the quantity of interest by constructing estimators that have smaller variance than does MCS. This technique was proposed in the 1950s \cite{Kahn_1953}. The major cause for the inefficiency in computing rare-event probabilities using MCS is  that most of the random samples generated do not belong to the extreme excursion region (or the region of interest). The basic idea of IS is to  use the information available about the rare event to generate samples that more frequently visit the region of interest. This is achieved by constructing an  IBD, which can be used to generate samples. If successful, unlike in the case of MCS, an appreciable fraction of these samples contribute to the probability estimate. When designing an IBD, the aim is for its probability mass to be concentrated in the region of interest. Based on this consideration, several techniques for constructing IBDs have been developed, such as variance scaling and mean shifting \cite{Bucklew_2013B}. A more detailed treatment of importance sampling and the relevant literature can be found in standard stochastic simulation textbooks \cite{Dunn_2011B, Asmussen_2007B}. One of the major challenges involved with importance sampling is the construction of an IBD that results in a low-variance estimator. We note that the approach may sometimes be inefficient in high-dimensional problems \cite{Katafygiotis_2008}. 
%%\cm{Vishwas, you need to reduce sections 1.1.3 to 1.1.7 by A LOT; you give a lot of details about methods that we will not use in the end or build on. So the referees would complain that it is besides the point. You may want to mention these alternatives, but primarily need to mention their LIMITATIONS compared to ours; particularily if you go to details in describing them. In the end you do not code and compare with ANY of them. If you describe them, why would you expect your approach to be better than theirs? }
%%\cv{I have summarized the methods without going into great detail in sections 1.1.3 and 1.1.4; I have also written briefly about their limitations and how our method goes some way in addressing them. I have also discussed other advantages of our methods (sec 1.2 and 1.3)}
\subsubsection{Nested subset methods}\label{sec:nested}The underlying idea of this class of methods is to consider a sequence of nested subsets of the probability space of interest (for example, starting with the entire space and shrinking to the target rare event) and use the notion of conditional probability to factorize the target event as a product of conditional events. Two main methods that fall into this class are subset simulation (SS) \cite{Au_2001} and splitting methods \cite{Kahn_1951}. In SS, the nested subsets are generated by choosing appropriate intermediate thresholds. Splitting methods are based on the idea of restarting the associated Markov process from certain system states in order to generate more occurrences of the rare event of interest. Several modifications have been proposed to both SS \cite{Ching_2005, Ching_2005_2, Katafygiotis_2005, Zuev_2012, Bect_2017} and splitting methods \cite{Botev_2012, Beck_2016}. Evaluating the conditional probabilities forms a major portion of the computational load. Computing the conditional probabilities for different nested subsets concurrently is nontrivial. Additionally, it is not clear how many model evaluations are required at the intermediate level sets in order to achieve a good probability estimate.
\subsubsection{Methods based on large deviation theory} Recent work by Dematteis et al.~used large deviations theory (LDT) to estimate the probabilities of rogue waves of a certain height \cite{Dematteis_2018}. The same authors used LDT to estimate probabilities of extreme events in dynamical systems with random components \cite{Dematteis_2019}. LDT is an efficient approach for estimating rare events  when the event space of interest is dominated by a few elements. 
%For examples demonstrated in \cite{Dematteis_2019}, the authors estimate the probability of rare events at a fixed time instant. Our problem statement is slightly different since we are interested in rare events in a time interval. 
The aforementioned papers solve an optimization problem to estimate the rare-event probability. In contrast, our approach uses a Bayesian inverse problem framework to determine an IBD, which will then be used to estimate the rare-event probability. In \S \ref{sec:ibd} we contrast the approach based on LDT with our approach.
\subsubsection{Multifidelity and surrogate-based methods} Multifidelity methods are used for estimating rare-event probabilities is situations when multiple evaluations of the forward model is prohibitively expensive. This approach leverages a hierarchy of low-cost reduced-order models, such as projection-based reduced-order models, data fit interpolation models, and support vector machines, to reduce the cost of constructing the IBD  \cite{Peherstorfer_2017}. The main idea behind the surrogate-based method is to start with a deterministic sample of the system and then construct a surrogate that approximates the system based on these samples \cite{Bucher_1990, Faravelli_1989, Wong_1985}. We remark that multifidelity and surrogate methods can be readily augmented with the framework developed in this paper to obtain additional computational savings.
\subsection{Overview of our methodology}Our methodology uses ideas from excursion probability theory to characterize the tails of the probability distribution associated with the event \cite{adler2010geometry}. Specifically, we use Rice's formula \cite{Rice_1944}, which was developed to compute the expected number of upcrossings for stochastic processes:
\begin{equation}\label{eqn:RicesFormula}
        \mathbb{E} \left\{N^{+}_u(0,T) \right \} = \displaystyle \int_{0}^{T} \, \int_{0}^{\infty}\, y \varphi_t(u,y)\, \mathrm{d}y\, \mathrm{d}t\, .
\end{equation}
%
%\cm{Later on we use $f$ for the dynamical system. We need to change this $f(t)$ to something else}
The left-hand side denotes the number of upcrossings of level $u$,  $y$ is the derivative of the stochastic process (in a mean square sense), and $\varphi_t(u,y)$ represents the joint probability distribution of the process $g(t)$ and its derivative $\frac{dg}{dt}$. Clearly the expression in the integral is analytically tractable only for special types of stochastic processes. Specifically, for Gaussian processes, this term can be resolved analytically. Moreover---the critical feature we will use here--- for smooth Gaussian processes $g(t)$, Rice's formula is a faster-than-exponentially-accurate approximation of the excursion probability. 
That is \cite[Equation (14.0.2)]{adler2009random}:
\begin{equation}\label{eq:RiceConvergence}
\left| \mathbb{P}\left\{ \sup_{t \in [0,T]} g(t) \geq u \right\} -  \mathbb{E} \left\{N^{+}_u(0,T) \right \}\right|
\leq  \mathcal{O}\left( e^{-\beta u^2}\right), 
\end{equation}
where $\beta > 0 $ is a parameter depending on the process $g(t)$ and interval $T$, but not on the target level $u$, and the asymptotics in the $\mathcal{O}()$ notation refers to $ u \rightarrow \infty$. If we use the number of upcrossings in \eqref{eqn:RicesFormula} as our estimate of the excursion probability, we can  interpret large values of $y \varphi_t(u,y)$ as defining the times and values of the process velocity for which the crossing is most likely to occur. This, in turn, is the key in efficiently determining the points in the input space that represent the highest contribution to the excursion probability. 

The setup in this article involves a nonlinear dynamical system that is excited by a Gaussian or a non-Gaussian initial state that results in a non-Gaussian stochastic process. 
%To address this problem, we (1) use Gaussian mixture models (GMM) to approximate the probability distribution of the initial state of the system (which we assume given). 
To address this problem, we  linearize the nonlinear dynamical system variation around the trajectories starting at the mean of the initial state. We thus obtain a Gaussian approximation to the system trajectory distribution. 
%We note that, if the input GMM distribution has sufficiently low variance components, the linear approximation will be accurate. 
Furthermore, we  use Rice's formula and solve a sequence of Bayesian inverse problems \eqref{eqn:RicesFormula} to determine the uncertainty sets in the input space that are most likely to cause the extreme events; these sets,  in turn, are used to construct the biasing distribution. The main advantages of our approach are the following:
\begin{itemize}
\item Constructing the biasing distribution is the most expensive component of the computational method. Since our method does not use nested subsets to evaluate the target probability, it is amenable to parallelization (see the discussion in \S \ref{sec:mcmc-based-sampling}).
%\item The computational framework presented in this work is agnostic to the probability distribution of the random variables (though its complexity may be influenced by it, through the number of GMM components required for approximating the input distribution to a given level of precision).
\item As we will demonstrate in \S \ref{sec:num_exp}, a moderate number of evaluations of the model ($\mathcal{O}(1000)$) are required in order to achieve acceptable levels of relative accuracy ($\mathcal{O}(10^{-2})$). The method can capture probabilities on the order of $\mathcal{O}(10^{-6})$ accurately.% and the performance gap between MCS and our method for lower probabilities is even more significant as demonstrated in \S \ref{sec:num_exp}.
\item Although we use nonlinear dynamical systems with random initial states as a basis to demonstrate our method, the algorithm can be seamlessly extended to stochastic dynamical systems with random parameters.
\item In applications that require repeated evaluations of the rare-event probability and where the distribution of the random parameter does not change significantly between these evaluations, the IBD can be reused to obtain accurate estimates of the rare-event probability. We demonstrate this for a small problem in \S \ref{sec:num_exp}.
\end{itemize}
\subsection{Limitations and possible extensions} One of the major limitations of our approach is that as the dimensionality of the random variable grows, the construction of the biasing distribution becomes expensive. Currently, this method is practical for problems where the size of the random variable is $\mathcal{O}(100)$. However, we aim to solve problems up to and beyond $\mathcal{O}(1000)$. The current approach requires that the random variable be normally distributed. For many practical problems, however, this might not be the case. In such scenarios we use the method of moments to approximate the non normal random variable by a Gaussian distribution. Another possible approach to handling non-Gaussian random parameters is to use a Gaussian mixture model (GMM) to approximate it. In such a scenario, challenges may arise regarding controlling the variance of the GMM components such that the errors due to linearization do not grow too much. An obvious extension to the current work is to develop a strong theoretical foundation that justifies the algorithm in this paper. Another potential research direction is related to constructing the likelihood function that is necessary for the MCMC step of the algorithm. Currently, we use ad hoc methods to choose the likelihood for the MCMC step (more details are in \S \ref{sec:mcmc-based-sampling}), and this approach.can be significantly improved by using a design of experiments approach (similar to \cite{Mohamad_2018}).
%of nonlinear dynamical systems, due to the linearization the accuracy estimates might be impacted if we fail to control the variance of the GMM components. Currently, we use the expectation maximization algorithm that is part of the statistical toolbox in MATLAB to construct the Gaussian mixture. However, it is desirable to develop a method where we have greater control on the eigenvalues of the covariance matrix associated with the GMM components. This can also aid us in addressing the aforementioned dimensionality issues. 

\section{The rare-event problem}\label{sec:problem} Consider an input-output system with $d \in \mathbb{N}$ inputs and $d' \in \mathbb{N}$ outputs, which is modeled by a continuous function $s: \mathbb{R}^d \rightarrow \mathbb{R}^{d'}$. 
%Let $(\Omega, \mathcal{F}, \mathbb{P})$ be a probability space with the sample space $\Omega$, the $\sigma-$algebra $\mathcal{F}$, and the probability measure $\mathbb{P}$. 
We represent the uncertainties in the input by a probability distribution with probability density function (PDF) $p$. Let $Z$ be a random $d$-dimensional vector. The inputs to the system are the components $\mathbf{z} = [z_1, \ldots, z_d]^{\top} \in \mathbb{R}^d$ of a realization of $Z$. The outputs $s(\mathbf{z}) \in \mathbb{R}^{d'}$ are the realization of the random variable $s(Z)$. We are interested in the failure probability of the system described by the model $s$. Let $g:\mathbb{R}^{d'} \rightarrow \mathbb{R}$ be the limit-state function, which we assume to be continuous. We say that the system fails for an input $\mathbf{z}$ if $g(s(\mathbf{z})) \geq \mathbf{U}$. This leads us to the failure domain $$\mathcal{F} = \left \{ \mathbf{z} \in \mathbb{R}^d : g(s(\mathbf{z})) \geq \mathbf{U}\right\}\,$$ and the indicator function of the failure domain 
\begin{align}
        \mathbb{I}_{\mathcal{F}} =
\begin{cases}
1\,, \qquad \mathbf{z} \in \mathcal{F}\\
0\, \qquad \mathbf{z} \not \in \mathcal{F}\,.
\end{cases}
\end{align}
We define the failure probability of the system as 
\begin{align}\label{eqn:expectationProbabilites}
        P_{\mathcal{F}} = \mathbb{E}_p \left \lbrack \mathbb{I}_{\mathcal{F}}(Z)\right \rbrack = \displaystyle \int_{\mathbb{R}^d} \, \mathbb{I}_{\mathcal{F}}(\mathbf{z}) p (\mathbf{z})\, \mathrm{d} \mathbf{z}\,. 
\end{align} 
The variance of $\mathbb{I}_\mathcal{F}(Z)$ with respect to the PDF $p$ is 
\begin{align}\label{eqn:variance}
        \textrm{var}_p \left \lbrack \mathbb{I}_{\mathcal{F}}(Z)\right \rbrack = \displaystyle \int_{\mathbb{R}^d} \,\left(\mathbb{I}_{\mathcal{F}} - \mathbb{E}_p \left \lbrack \mathbb{I}_{\mathcal{F}}(Z) \right \rbrack \right)^2 p(\mathbf{z})\, \mathrm{d}\mathbf{z} = P_{\mathcal{F}} - P_{\mathcal{F}}^2\,.
\end{align}
The MCS method is often used to estimate expectation such as \eqref{eqn:expectationProbabilites}. It draws $M \in \mathbb{N}$ independent and identically distributed (i.i.d.) samples $\mathbf{z}_1, \ldots, \mathbf{z}_M \in \mathbb{R}^d$ from the distribution of $Z$, that is, $M$ realizations of the random variable $Z$, and computes the Monte Carlo estimate
\begin{align}\label{eqn:Expectation_MC}
P_{\mathcal{F}}^{MC}(\mathbf{z}_1, \ldots, \mathbf{z}_M) = \frac{1}{M} \sum_{i=1}^{M}\mathbb{I}_{\mathcal{F}} (\mathbf{z}_i)\,.
\end{align}
Note that we distinguish between the estimate $P_{\mathcal{F}}^{MC}(\mathbf{z}_1, \ldots, \mathbf{z}_n)$, which is a scalar value, and the Monte Carlo estimator $P_{\mathcal{F}}^{MC}(Z)$, which is a random variable. The relative root mean square error (RMSE) of $P_{\mathcal{F}}^{MC}(Z)$ is 
\begin{align}\label{eqn:mcerror}
e(P_{\mathcal{F}}^{MC}) \approx \frac{1}{\mathbb{E}_p \left \lbrack \mathbb{I}_{\mathcal{F}}(Z)\right \rbrack}\sqrt{\frac{\textrm{var}_p \left \lbrack \mathbb{I}_{\mathcal{F}(Z)} \right \rbrack}{M}} = 
\sqrt{\frac{P_{\mathcal{F}}-P_{\mathcal{F}}^2}{P^2_{\mathcal{F}}M}} \approx \sqrt{\frac{1}{P_{\mathcal{F}} M}}
\,,
\end{align}
when $P_{\mathcal{F}} \ll 1$.
Hence, for a threshold parameter $0 < \epsilon \in \mathbb{R}$, a relative error  $e(P_{\mathcal{F}}^{MC}) \leq \epsilon$ is achieved with 
\begin{align}
        M = \left\lceil \frac{1}{P_{\mathcal{F}} \epsilon^2}\right\rceil
\end{align}
samples, where $ \lceil \rceil$ denotes the ceil function. If $P_{\mathcal{F}} \ll 1$, then $M$ is large. Hence, estimating small probability events with MC is difficult.
\subsection{Importance sampling}\label{sec:is}
Variance reduction methods aim to reduce the RMSE \eqref{eqn:mcerror} by changing  \eqref{eqn:expectationProbabilites} to a new integral with the same value but with an integrand and/or a distribution that, combined, result in a lower variance than the original function $\mathbb{I}_{\mathcal{F}}$ has with respect to the distribution $p$. Importance sampling is one such variance reduction method that has been used successfully to estimate failure probabilities\cite{Srinivasan_2013B, Robert_2005B}. Importance sampling introduces a random vector $Z': \Omega \rightarrow \mathbb{R}^d$ with PDF $q$, which is used as the biasing distribution. A realization of $Z'$ is denoted by $\mathbf{z}' = \left \lbrack z_1', \ldots, z_d'\right \rbrack^{\top} \in \mathbb{R}^d$. In the following, the distribution of $Z$ is the nominal distribution, and the corresponding PDF $p$ is the nominal PDF. The distribution of $Z'$ is the biasing distribution, and $q$ is the biasing PDF. The biasing PDF $q$ is constructed such that the $\textrm{supp}(p)\subseteq \textrm{supp}(q)$, where $$\textrm{supp}(p) = \{\mathbf{z}\in \mathbb{R}^d : p(\mathbf{z}) > 0\}\,,$$ denotes the support of the PDF $p$. Let $w : \textrm{supp}(p) \rightarrow \mathbb{R}$ be the weight function $w(\mathbf{z}') 
\stackrel{\Delta}{=} \displaystyle \frac{p(\mathbf{z}')}{q(\mathbf{z}')}$. The weight $w(\mathbf{z}')$ is the importance weight of a realization $\mathbf{z}' = Z'(\omega)$. Because $\textrm{supp}(p) \subseteq \textrm{supp}(q)$ holds, the failure probability \eqref{eqn:expectationProbabilites} equals the expectation of the random variable $\mathbb{I}_{\mathcal{F}}(Z')$ weighted with the random variable $w(Z')$. That is, we have 
\begin{align}\label{eqn:importance_expectation}
        P_{\mathcal{F}} = \mathbb{E}_q \left \lbrack \mathbb{I}_{\mathcal{F}}(Z')w(Z')\right \rbrack\,.
\end{align}
The expectation \eqref{eqn:importance_expectation} is approximated with the Monte Carlo method, with samples $\mathbf{z}_1', \ldots, \mathbf{z}_M'$ drawn from the biasing distribution. Thus the importance sampling estimate of $P_{\mathcal{F}}$ with samples $\mathbf{z}_1', \ldots, \mathbf{z}_M'$ is
\begin{align}
P_{\mathcal{F}}^{\textrm{IS}}(\mathbf{z}_1', \ldots, \mathbf{z}_M') = \frac{1}{M} \sum_{i=1}^{M} \mathbb{I}_{\mathcal{F}} (\mathbf{z}_i')w(\mathbf{z}_i')\,.
\end{align}
Therefore, the Monte Carlo method with importance sampling consists of two steps. In step one, the biasing distribution is generated. In step two, the importance sampling estimator $P_{\mathcal{F}}^{\textrm{IS}}(Z')$ is an unbiased estimator of $P_{\mathcal{F}}$ because $\textrm{supp}(p) \subseteq \textrm{supp}(q)$ \cite{Robert_2005B}.

If $\mathbb{E}_q \left \lbrack \mathbb{I}_{\mathcal{F}}(Z')^2 w(Z')^2 \right \rbrack < \infty$, then the relative RMSE of the importance sampling estimator is 
\begin{align}\label{eqn:is_variance}
e(P_{\mathcal{F}}^{\textrm{IS}}) = \frac{1}{P_{\mathcal{F}}}\sqrt{\frac{\textrm{var}_q\left \lbrack \mathbb{I}_{\mathcal{F}}(Z') w(Z') \right \rbrack}{M}}\,.
\end{align}
If the variance $\textrm{var}_q\left \lbrack \mathbb{I}_{\mathcal{F}}(Z') w(Z') \right \rbrack$ is smaller than $\textrm{var}_p \left \lbrack \mathbb{I}_{\mathcal{F}}(Z) \right \rbrack, $ then the relative RMSE of the importance sampling estimator is smaller than the relative RMSE of the Monte Carlo estimator for the same number of samples $M$.

\section{Construction of IBD via Bayesian inference} \label{sec:mcmc-based-sampling}
Consider the following dynamical system,
\begin{align}\label{eqn:dynamics}
        \x' &=  f(t,\x)\,, \quad t  =[0, T] \\
        \x(0) &= \x_0\,, \quad \x_0 \sim p\,, \quad \x \in \Omega\,, \nonumber
\end{align}
% \[
%         \x' =
%         \begin{cases}\label{eqn:dynamics}
%                 f(t,\x)\,, \quad t  =[0, T] \\
%                 \x(0) = \x_0\,, \quad \x_0 \sim p\,, \quad \x \in \Omega\,, 
%         \end{cases}
% \]
% \begin{cases}\label{eqn:dynamics}
%         \x' = f(t,\x)\,, & t  \\
%         \quad \x(t_0) = \x_0\,, \quad \x_0 \sim p\,, \x \in \mathbb{R}^d\,, 
% \end{cases}
%
where the initial state of the system $\x_0$ is uncertain and has a probability distribution $p$, with $\mu$ being the corresponding probability measure. The problem of interest to us is to estimate the probability that $ \mathbf{c}^{\top} \mathbf{x}(t)$ exceeds the level $u$ for $t\in \lbrack 0, T \rbrack$. That is, we seek to estimate the {\textit{excursion probability}}
\begin{align}\label{eqn:problem_of_interest}
P_T(u) := \mathbb{P}\left( \sup_{0 \leq t \leq T} \mathbf{c}^{\top}\mathbf{x}(t, \x_0) \geq u\,, ~~t\in \lbrack 0, T \rbrack \right)\,,
\end{align}
where $\x(t, \x_0)$ represents the solution of the dynamical system \eqref{eqn:dynamics} for a given initial condition $\x_0$.
% In other words, we are interested in estimating the probability that $\mathbf{c}^{\top} \mathbf{x}(t)$ exceeds the level $u$ for $t$ in the interval $\lbrack 0, T \rbrack$. 
% Let 
% \begin{align}\label{eqn:solution_initial_state}
%         \mathbf{x} = g(t, \x_0)
% \end{align}
% represent the solution to the dynamical system $\displaystyle \x' = f(t,\x)$ in terms of the initial state of the system. Also 
Let $\Omega(u) \subset \Omega$ represent the set of all initial conditions for which the solution of the dynamical system exceeds the excursion level $u$. That is,
\begin{align}
\Omega(u) := \left \{\x_0 : \sup_{0 \leq t \leq T}\mathbf{c}^{\top}\x(t,\x_0) \geq u  \right \}\,.
\end{align}
Notice that $\x(t, \x_0)$ depends on $\x_0$ implicitly through the solution of the dynamical system.
%  defined by $\left \{\forall \x_0 \in \mathcal{B}_u\,, \mathbf{c}^{\top} g(t,\x_0) \geq  u, \text{ for some } t \in \lbrack0, T \rbrack \right \}$. 
Since we can write 
\begin{align}\label{eqn:probability_conditional}
   %\displaystyle P_u = \frac{P(\x_0 \in \mathcal{B}_u \cap p)}{P(\x_0 \in p)}\,.
   P_{T}(u) = \mu(\Omega(u))\,,
\end{align}
estimating $P_T(u)$ is related to determining $\Omega(u)$. In general, one cannot determine $\Omega(u)$ analytically. We use Rice's formula, \eqref{eqn:RicesFormula}, to gain insight about $\Omega(u)$, which in turn will be used to construct an approximation to $\Omega(u)$. Let us revisit Rice's formula:
\begin{equation}\label{eqn:RicesFormula_2}
        \mathbb{E} \left\{N^{+}_u(0,T) \right \} = \displaystyle \int_{0}^{T} \, \int_{0}^{\infty}\, y \varphi_t(u,y)\, \mathrm{d}y\, \mathrm{d}t\, .
\end{equation}
{Recall that $\varphi_t(u,y)\,$ represents the joint probability density of $\mathbf{c}^{\top}\x$ and its derivative $\mathbf{c}^{\top}\x'$ for an excursion level $u$. The right-hand side of equation \eqref{eqn:RicesFormula_2} integrates the joint density over all values of derivatives and times at which there is an excursion. The key insight for our method is that values of the time $t$ and slope $y$ at which $y \varphi_t(u,y)\,$ is large contribute the most to this integral. We use this idea 
%suggests likely regions of excursion %Evaluating $y \varphi_t(u,y)\,$ for different values of $y$ and $t$ gives us an intuition about the times and the corresponding slopes at which the stochastic process potentially exceeds the value $u$ 
to construct an approximation $\widehat{\Omega}(u)$ to $\Omega(u)$.} 

Using \eqref{eqn:RicesFormula_2}, we can interpret $y \varphi_t(u,y)\,$ as an unnormalized PDF and thus sample from it to compute $\mathbb{E} \left\{N^{+}_u(0,T) \right \}$ using Monte Carlo approximation.
By sampling from the unnormalized distribution $y \varphi_t(u,y)\,$, we obtain a slope-time pair, $(y_i, t_i)$ at which the sample paths of the stochastic process exceed the excursion level $u$. Consider the forward map $\mathcal{G}: \mathbb{R}^{d\times 1} \rightarrow \mathbb{R}^2$, which evaluates the vector $\displaystyle \begin{bmatrix}\mathbf{c}^{\top}\x(t) \\ \mathbf{c}^{\top}\x'(t) \end{bmatrix}$ based on the dynamics \eqref{eqn:dynamics}, given an initial state $\x_0$ and a time $t$. We call $\x_i$ a preimage of a sample $(y_i, t_i)$ if 
        \begin{align}
        \mathcal{G}(\x_i, t_i) = \displaystyle \begin{bmatrix} u \\ y_i \end{bmatrix} \,.   
        \end{align}
Note that the problem of finding a preimage of a sample $(y_i, t_i)$ is ill-posed. Multiple $\x_i$'s  map to $ \displaystyle \begin{bmatrix} u \\ y_i \end{bmatrix}$ at time $t_i$ via the operator $\mathcal{G}$. Therefore, we define the set 
\begin{align}
        X_i := \left \{\x_i \in \Omega : \mathcal{G}(\x_i, t_i) = \displaystyle \begin{bmatrix} u \\ y_i \end{bmatrix}  \right \}\,,
\end{align}
and we construct our approximation $\widehat{\Omega}(u)$ as
\begin{align}\label{eqn:approximate_omega_u}
        \widehat{\Omega}(u) := \bigcup_{i=1}^N X_i\,.
\end{align}
Our intuition is that $\widehat{\Omega}(u)$ approximates $\Omega(u)$ better as we increase $N$. The underlying computational framework to approximate $\widehat{\Omega}(u)$ consists of the following stages:
\begin{itemize}
        \item Draw samples from unnormalized $y \varphi_t(u,y)\,$
        \item Find the preimages of these samples to approximate $\Omega(u)$.
\end{itemize}
We use MCMC to draw samples from unnormalized $y \varphi_t(u,y)\,$. We note that irrespective of the size of the dynamical system, $y \varphi_t(u,y)\,$ represents an unnormalized density in two dimensions; hence, using MCMC is an effective means to draw samples from it. Drawing samples from $y \varphi_t(u,y)\,$ requires evaluating it repeatedly, and in the following section we discuss the means to do so.
\subsection{Evaluating $y \varphi_t(u,y)\,$}\label{sec:evaluatevarphi}  %The joint probability density of the stochastic process and its derivative, $y \varphi_t(u,y)$, needs to be evaluated repeatedly in our method. 
%To construct the likelihood PDF, we must construct $\y_i$. In order to construct $\y_i$, we need to sample from the unnormalized density $y \varphi_t(u,y)$. 
In this section, we describe the process of evaluating $y \varphi_t(u,y)$ given $y$, $t$, and $u$. We note that $y \varphi_t(u,y)\,$, can be evaluated analytically only for special cases. Specifically, when $\varphi_t(u,y)$ is a Gaussian process, the joint density function $y \varphi_t(u,y)\,$ is analytically computable. Consider the dynamical system described by \eqref{eqn:dynamics}. When $p$ is Gaussian and $f$ is linear, we have
\begin{align}\label{eqn:linearGaussian}
        \x' = A\, \x(t) + b\,, \quad \x(t_0) = \x_0\,, \quad \x_0 \sim \mathcal{N}(\overline{\x}_0, \Sigma)\,.
\end{align}
Assuming $A$ is invertible, we can write $\x(t)$  as
\begin{align}\label{eqn:closedform}
\x(t) = \exp(A(t-t_0)) \, \x_0  - \left(I  - \exp(A(t-t_0))\right)A^{-1} b\,,
\end{align}
where $I$ represents an identity matrix of the appropriate size. Given that $\x_0$ is normally distributed, it follows that $\x(t)$ is a Gaussian process:
\begin{align}\label{eqn:gpx}
        &\x(t) \sim \mathcal{GP} \left(\overline{\x}, {\rm cov}_{\x} \right ) \,, \text{ where }\\
     &\nonumber     \overline{\x} = \exp(A(t-t_0)) \overline{\x}_0 - \left(I  - \exp(A(t-t_0))\right)A^{-1} b \, \text{ and } \\
     &\nonumber   {\rm cov}_{\x} = \exp(A(t-t_0)) \Sigma \left(\exp(A(t-t_0))\right)^\top\,.
        % \x(t) \sim \mathcal{GP}\left(\exp(A(t-t_0)) \mu - \left(I  - \exp(A(t-t_0))\right)A^{-1} b, \exp(A(t-t_0)) \Sigma \left(\exp(A(t-t_0))\right)^\top\right)
\end{align}
The joint PDF of a stochastic process and its derivative, 
 $\varphi$, the joint PDF of $\mathbf{c}^{\top}\mathbf{x}(t)$ and $\mathbf{c}^{\top}\mathbf{x}'(t)$, is given by \cite[equation 9.1]{3569}
\begin{align}\label{eqn:gpvarphi}
        \begin{bmatrix}
                \mathbf{c}^{\top}\x \\
                \mathbf{c}^{\top}\x'
        \end{bmatrix}
        \sim \mathcal{GP}\left(\overline{\x}^{\varphi}, \begin{bmatrix}
                \mathbf{c}^{\top}\Phi\mathbf{c} & \mathbf{c}^{\top}\Phi A^{\top}\mathbf{c} \\
                                            \mathbf{c}^{\top}A \Phi^{\top}\mathbf{c}& \mathbf{c}^{\top}A\Phi A^{\top}\mathbf{c}
                                           \end{bmatrix}
        \right)\,,  
\end{align}
where $$\overline{\x}^{\varphi} := \begin{bmatrix} \mathbf{c}^{\top} \overline{\x} \\ \mathbf{c}^{\top} (A\overline{\x} + b)\end{bmatrix}$$ and $$\Phi := \exp(A(t-t_0)) \Sigma \left(\exp(A(t-t_0))\right)^\top\,.$$
We can now evaluate $y \varphi_t(u,y)\,$ for arbitrary values of $u_i$, $y_i$, and $t_i$ as 
\begin{align}\label{eqn:compute_varphi}
y_i \varphi_{t_i}(u_i,y_i) = \frac{y_i}{2 \pi \mid \Upsilon \mid}\exp\left(-\frac{1}{2} \left \| \begin{bmatrix} u_i \\ y_i \end{bmatrix} - \overline{\x}^{\varphi}\right \|^2_{\Upsilon^{-1}}\right)\,,
\end{align}
where $\Upsilon :=  \begin{bmatrix}
        \mathbf{c}^{\top}\Phi\mathbf{c} & \mathbf{c}^{\top}\Phi A^{\top}\mathbf{c} \\
                                    \mathbf{c}^{\top}A \Phi^{\top}\mathbf{c}& \mathbf{c}^{\top}A\Phi A^{\top}\mathbf{c}
                                   \end{bmatrix}$ and $\mid \Upsilon \mid$ denotes the determinant of $\Upsilon$. Note that the right-hand side in \eqref{eqn:compute_varphi} is dependent on $t_i$ via $\Upsilon$.
\subsubsection{Notes for nonlinear $f$}\label{sec:nonlinearf} When $f$ is nonlinear,  one cannot compute $y \varphi_t(u,y)\,$ analytically---a key ingredient for our computational procedure. We approximate the nonlinear dynamics by linearizing $f$ around the mean of the initial distribution. Assuming that the initial state of the system is normally distributed as described by equation \eqref{eqn:linearGaussian}, linearizing around the mean of the initial state gives 
        \begin{align}\label{eqn:linearizedsystem}
               \x' \approx \mathbf{F} \cdot (\x - \overline{\x}_0) +f(\overline{\x}_0, 0) \,,
               \end{align}
       where $\mathbf{F}$ represents the Jacobian of $f$ at $t=0$, $\x=\overline{\x}_0$. This reduces the nonlinear dynamical system to a form that is similar to equation \eqref{eqn:linearGaussian}. Thus, we can now use equations \eqref{eqn:gpx}, \eqref{eqn:gpvarphi}, and \eqref{eqn:compute_varphi} to approximate  $y \varphi_t(u,y)\,$ for nonlinear $f$.

We now describe a systematic computational framework to determine $X_i$ for a given sample $(y_i, t_i)$. This allows us to determine the elements of set $\widehat{\Omega}(u)$. 
%\subsection{Constructing $X_i$ for a given sample $(y_i, t_i)$}\label{sec:approxB}
\subsection{Determining the preimages for a given sample}\label{sec:approxB}
A sample from the unnormalized joint distribution $y \varphi_t(u,y)\,$ gives a slope, $y_i$, and time, $t_i$, at which the stochastic process exceeds the level $u$. Hence $\displaystyle \begin{bmatrix}\mathbf{c}^{\top}\x(t_i) \\ \mathbf{c}^{\top}\x'(t_i) \end{bmatrix} =  \begin{bmatrix} u \\ y_i \end{bmatrix} $. Constructing $X_i$ requires finding all the preimages $\mathcal{G}^{-1} \left(\begin{bmatrix} u \\ y_i \end{bmatrix}\right)  \subset \Omega$. This amounts to finding all the solutions of the following equation, 
\begin{align}\label{eqn:nonlinear}
        \mathcal{G}(\x, t_i) = \y_i\,,
\end{align}
where $ \y_i = \begin{bmatrix} u \\ y_i \end{bmatrix}$. Another formulation of the problem \eqref{eqn:nonlinear} is
        \begin{equation}\label{eqn:ip_sample}
                \begin{aligned}
               \x_i := & \underset{\x}{\textrm{  arg min}}
                & & \frac{1}{2} \| \y_i - \mathcal{G}(\x, t_i) \|_2^2 \,.
                \end{aligned}
        \end{equation}
        Since $\mathcal{G}$ is a mapping from $\mathbb{R}^{d \times 1}$ to $\mathbb{R}^2$, problem \eqref{eqn:ip_sample} is an ill-posed and underdetermined inverse problem. To address the ill-posedness, we use the Bayesian formulation of the inverse problem by placing a prior on $\x_{i}$ and identifying the term $\| \y_i - \mathcal{G}(\x, t_i) \|_2^2 $ as a negative log-likelihood. Suppose, in the process of finding preimages $\mathcal{G}^{-1} \left(\begin{bmatrix} u \\ y_i \end{bmatrix}\right)$, we encounter elements that map to a value higher than $u$. These should not be discarded because these elements still cause an excursion and hence are elements of the set $\Omega(u)$. The Bayesian treatment allows for such flexibility because of the covariance associated with the log-likelihood term. We note, however, that the nonlinear equation \eqref{eqn:nonlinear} does not allow this flexibility. 
        
        In equation \eqref{eqn:dynamics}, we stated that $\x_0$ has a probability distribution $p$ and that we use $p$ as a prior PDF for  $\x_i$. 
        %Assuming $p$ is Gaussian with mean $\mathbf{0}$ and covariance $\Sigma$, we have
        % \begin{align}\label{eqn:prior}
        %         \pi_{\rm pr}(\x_i) \propto \exp\left(-\frac{1}{2} \left \lVert \x_i \right \rVert^2_{\Sigma^{-1}}\right)\,.
        % \end{align}
        \begin{align}
                \pi_{\rm pr}(\x_i) \propto p
        \end{align}
        Treating $\y_i$ as a random variable with covariance $\Gamma_i$, we can write the following likelihood:
        \begin{align}\label{eqn:likelihood}
                \pi_{\rm like} (\y_i \mid \x_i ) \propto \exp \left (-\frac{1}{2} \left \lVert{\y}_i - \mathcal{G}( \x_i, t_i) \right \rVert^2_{\Gamma^{-1}_i}\right)\,.
        \end{align}
        Using Bayes' rule, we can write the posterior PDF of $\x_i$ as 
\begin{align} \label{eqn:posterior}
        \pi_{\rm post}^i :=   \pi_{\rm post}( \x_i \mid \y_{i}) \propto p \, \pi_{\rm like} (\y_{i} \mid \x_i )\,,
\end{align}
which is
% \begin{align}\label{eqn:posterior2}
%       \pi_{\rm post}^i( \x_i \mid \y_{i}) \propto \exp\left (-\frac{1}{2} \left \lVert \x_i \right \rVert^2_{\Sigma^{-1}}\right) \, \exp \left (-\frac{1}{2} \left \lVert\overline{\y}_i - \mathcal{G}( \x_i, t_i) \right \rVert^2_{\Gamma^{-1}_i}\right)\,.
%         %\pi_{\rm post}( \x^i_{\mathcal{B}} \mid \y^{i}_{\rm obs}) \propto \exp\left (-\frac{1}{2} \left \lVert \x^i_{\mathcal{B}} \right \rVert^2_{\Sigma^{-1}} -\frac{1}{2} \left \lVert \mathcal{G}(\x^i_{\mathcal{B}}) - \y^{i}_{\rm obs} \right \rVert^2_{{\Gamma_{\rm obs^i}^{-1}}}  \right )\,.
% \end{align}
\begin{align}\label{eqn:posterior2}
              \pi_{\rm post}^i( \x_i \mid \y_{i}) \propto p \, \exp \left (-\frac{1}{2} \left \lVert{\y}_i - \mathcal{G}( \x_i, t_i) \right \rVert^2_{\Gamma^{-1}_i}\right)\,.
                %\pi_{\rm post}( \x^i_{\mathcal{B}} \mid \y^{i}_{\rm obs}) \propto \exp\left (-\frac{1}{2} \left \lVert \x^i_{\mathcal{B}} \right \rVert^2_{\Sigma^{-1}} -\frac{1}{2} \left \lVert \mathcal{G}(\x^i_{\mathcal{B}}) - \y^{i}_{\rm obs} \right \rVert^2_{{\Gamma_{\rm obs^i}^{-1}}}  \right )\,.
        \end{align}
When $p$ is Gaussian, the kernel of the likelihood distribution can be represented in closed form. Hence the Bayesian inverse problem in \eqref{eqn:posterior2} can be solved either by finding a maximum a posteriori point (MAP) and using the Laplace approximation around the MAP to describe the uncertainty around the solution or by drawing samples from the approximate posterior distribution using MCMC. In scenarios when $p$ is non-Gaussian, however, the challenges are twofold:
\begin{itemize}
        \item The kernel of the posterior cannot be represented in a closed analytical form.
        \item We cannot evaluate $y\varphi_t(u, y)$ analytically---which is central to our method.
\end{itemize} 
We tackle non-Gaussianity by using the method of moments to approximate $p$ by a Gaussian distribution. Using Gaussian mixtures might lead to a better approximation to $p$ than using the method of moments, and we can reuse the technique here about the center of each component of the mixture. However, this approach suffers from the curse of dimensionality when $p$ represents a PDF in large dimensions, and hence we do not pursue using Gaussian mixtures to approximate $p$ in this paper, with the expectation that our approach will create an acceptable IBD (which need not be exact). Assuming $p$ is Gaussian or can be approximated by the method of moments,  we can write $p$ as $\mathcal{N}(\overline{\x}, \Sigma)$. 
\begin{align}\label{eqn:posterior2_approximate}
      \pi_{\rm post}^i( \x_i \mid \y_{i}) \propto \exp\left (-\frac{1}{2} \left \lVert \x_i - \overline{\x} \right \rVert^2_{\Sigma^{-1}}\right) \, \exp \left (-\frac{1}{2} \left \lVert{\y}_i - \mathcal{G}( \x_i, t_i) \right \rVert^2_{\Gamma^{-1}_i}\right)\,.
        %\pi_{\rm post}( \x^i_{\mathcal{B}} \mid \y^{i}_{\rm obs}) \propto \exp\left (-\frac{1}{2} \left \lVert \x^i_{\mathcal{B}} \right \rVert^2_{\Sigma^{-1}} -\frac{1}{2} \left \lVert \mathcal{G}(\x^i_{\mathcal{B}}) - \y^{i}_{\rm obs} \right \rVert^2_{{\Gamma_{\rm obs^i}^{-1}}}  \right )\,.
\end{align}
The covariance information $\Gamma_i$ is necessary in order to evaluate the posterior PDF given $\x$. We discuss the choice of the covariance $\Gamma_i$ in the next subsection. %\S \ref{sec:covarianceGamma}. 
% To construct $\y_i$ we need to sample from the unnormalized density $y \varphi_t(u,y)\,$. This in turn requires evaluating $y \varphi_t(u,y)\,$ repeatedly. Below, we describe the procedure to evaluate $y \varphi_t(u,y)\,$.
% %
\subsection{Choice of covariance}\label{sec:covarianceGamma}
\begin{figure}
        \begin{center}
                {\includegraphics[width=0.55\linewidth]{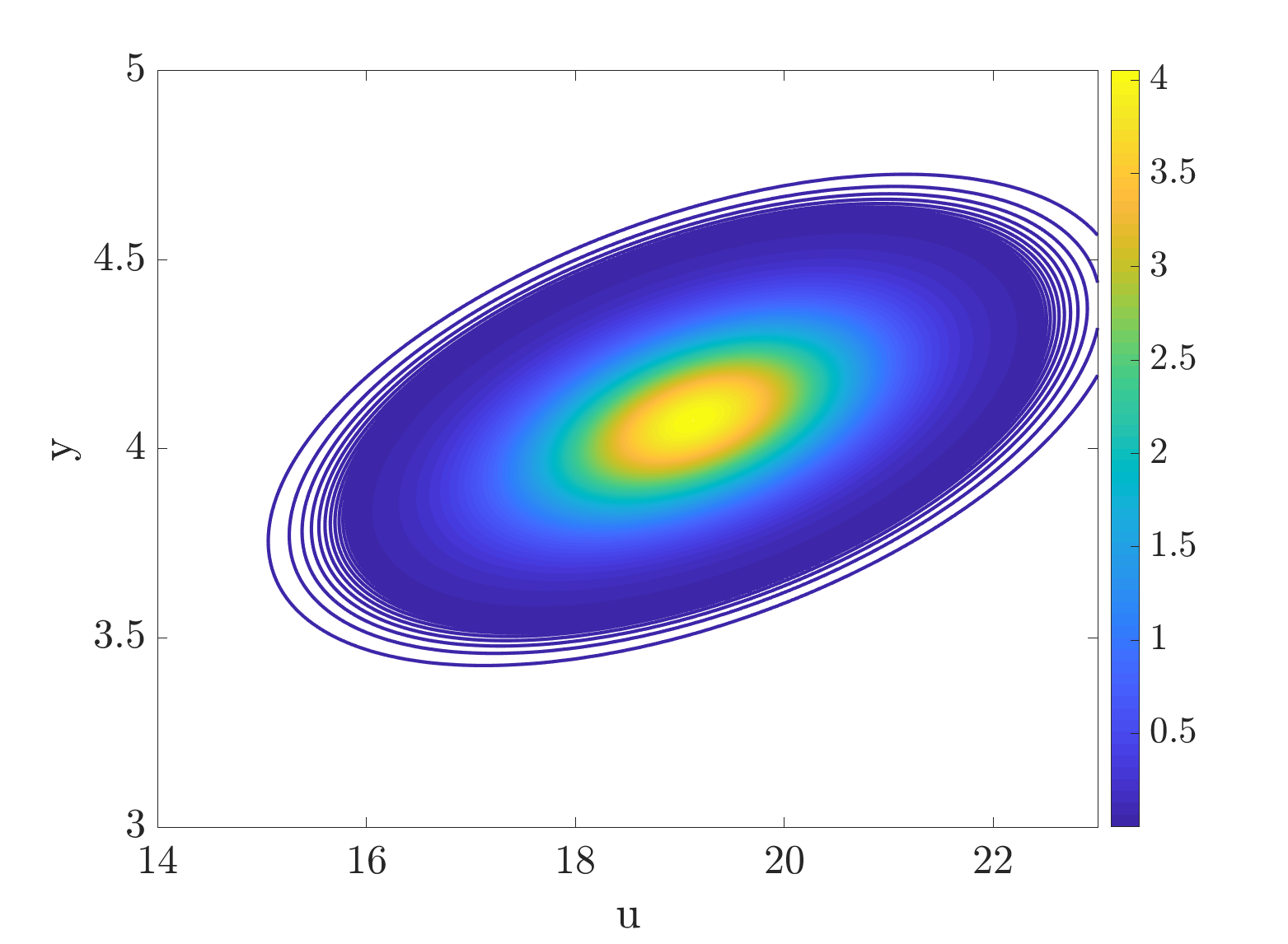}}
        \caption{Contours of $ y \varphi_t(u,y)$ evaluated for different values of $y$ and $u$ at a fixed $t_i$ for the Lotka-Volterra system.}
        \label{fig:distributionyphit}
        \end{center}
\end{figure}

For our specific problem, defining $\Gamma_i$ is an important step in solving the Bayesian inverse problem \eqref{eqn:posterior2_approximate}.
% We assume $\Gamma_i$ as diagonal. This assumption could lead to inefficiencies in our computational framework. However, unless we have problem specific insights about the correlation between $\mathbf{c}^{\top}\x(t_i)$ and $\mathbf{c}^{\top}\x'(t_i)$ it is difficult to make any statements about the off-diagonal elements of $\Gamma_i$. We also note that if such information is available, one can easily incorporate this in the covariance matrix. 
%We note that our choice of $\Gamma_i$ has a significant effect on the posterior distribution. 
%Suppose variance in the component $\mathbf{c}^{\top}\x(t_i)$ component of observation $\y_i^{\rm obs}$ is too small, then the samples in
We use the value of $y \varphi_t(u,y)$ as a guide to choose the covariance of the likelihood term in \eqref{eqn:posterior2_approximate}. Recall that a sample from unnormalized distribution $y \varphi_t(u,y)$ gives us a $(y_i, t_i)$ pair. To choose the covariance of $\y_i$, we look at the unnormalized distribution of $y \varphi_t(u,y)$ at time $t_i$.~That is, we model the joint distribution of $(u, y \mid t_i)$ based on the values of $y \varphi_t(u,y)$ evaluated at $t_i$. Specifically we evaluate $y \varphi_t(u,y)\mid_{t_i}$ for $[u, y] \in \lbrack u- \varepsilon_1, u + \varepsilon_1 \rbrack \times \lbrack y_i - \varepsilon_2, y_i + \varepsilon_2 \rbrack$. These values give us a range of slopes at time $t_i$ for which the state is close to the excursion level. We then use the values of $y \varphi_t(u,y)\mid_{t_i}$ to construct a Laplace approximation to obtain an approximation for the joint distribution of  $(u, y \mid t_i)$. This gives us an approximate covariance ($\Gamma_i$) for the likelihood PDF. This is illustrated for the Lotka-Volterra system in \Cref{fig:distributionyphit}. We evaluate $y \varphi_t(u,y)\mid_{t_i}$ at $t_i = 2$ for a range of values of $u$ and $y_i$ and fit a two-dimensional Gaussian distribution to approximate the covariance for the likelihood.
\section{Solution to the Bayesian inverse problem}\label{sec:BayesianInverseProblem}
In \S \ref{sec:mcmc-based-sampling} we formulated the process of approximating the $\Omega(u)$ as solving a sequence of Bayesian inverse problems. We also defined the necessary ingredients to define a Bayesian inverse problem---the prior and the likelihood function. Ideally, the solution to the Bayesian inverse problem in \eqref{eqn:posterior2_approximate} should yield the posterior distribution $\pi_{\rm post}^i( \x_i \mid \y_{i})$. Except under special circumstances, however, one cannot obtain a closed-form expression for the posterior distribution $\pi_{\rm post}^i( \x_i \mid \y_{i})$ \eqref{eqn:posterior2_approximate}. 
%(for example, when the forward map $\mathcal{G}$ is nonlinear a closed form expression for the posterior distribution is usually unknown). 
Let $\x_i^{\rm MAP}$ denote the maximum a posteriori point (MAP point), that is, the point that maximizes the posterior PDF (equation \eqref{eqn:posterior2_approximate}).
A standard approach to solving the Bayesian inverse problem \eqref{eqn:posterior2} is to first find $\x_i^{\rm MAP}$  and then approximate the forward map $\mathcal{G}$ by its linearization around $\x_i^{\rm MAP}$. This results in a Gaussian approximation to the posterior distribution $\pi_{\rm post}^i( \x_i \mid \y_{i}^{\rm obs})$, which is known as the Laplace approximation \cite{Tanbui_2013, Petra_2014}. 

Alternatively, one can use MCMC methods to sample from the posterior PDF. In the following paragraphs, we describe both these approaches for solving the Bayesian inverse problem.
\subsection{Laplace approximation at the MAP point} The problem of finding the point at which the posterior PDF \eqref{eqn:posterior2_approximate} is maximized can be formulated as a deterministic inverse problem. The negative log-likelihood is treated as the data misfit term, and the negative log prior is used as a regularizer to avoid overfitting. The resulting inverse problem can be written as
\begin{equation}\label{eqn:ip_laplace}
        \begin{aligned}
       \x^{\rm MAP}_i := & \underset{\x}{\textrm{  arg min}}
        & & \underbrace{\frac{1}{2} \| {\y}_i - \mathcal{G}(\x, t_i) \|_{\Gamma_i^{-1}}^2}_{\textrm{data misfit}} + \displaystyle \underbrace{\frac{\tau}{2}  \left \lVert \x  - \overline{\x}\right \rVert^2_{\Sigma^{-1}}}_{\textrm{regularization}}  \,,
        \end{aligned}
\end{equation}
where $\tau$ is the regularization parameter. The solution for the optimization problem in \eqref{eqn:ip_laplace} is the MAP point for the Bayesian inverse problem in \eqref{eqn:posterior2_approximate}. To solve the minimization problem in \eqref{eqn:ip_laplace}, we use gradient-based optimization methods (for example, L-BFGS);  the necessary gradient information can be evaluated by using adjoints. In Appendix \ref{sec:GradientHessian} we describe the computational procedure to evaluate the gradient information.

Assuming the forward map $\mathcal{G}(\cdot, \cdot)$ to be Fr\'echet differentiable, we can approximately express an observation $\y_i$ as 
\begin{align}
        \y_i \approx \mathcal{G}(\x_i^{\rm MAP}, t_i) + \displaystyle \frac{\partial \mathcal{G}}{\partial \x}(\x - \x_i^{\rm MAP}, t_i) + {\eta}\,,
\end{align}
where $\eta \sim \mathcal{N}(0, \Gamma_i)$ and $\displaystyle \frac{\partial \mathcal{G}}{\partial \x}$ is the Fr\'echet derivative of $\mathcal{G}$ evaluated at $(\x_i^{\rm MAP}, t_i)$. Hence the Laplace approximation of the posterior $\pi_i^{\rm post}(\x_i \mid \y_i)$ can be written as 
\begin{align}\label{eqn:laplace_approximation}
        \pi_i^{\rm post}(\x_i \mid \y_i) \sim \mathcal{N} (\x_i^{\rm MAP}, \Gamma_i^{\rm post})\,,    
\end{align}
where $\displaystyle \Gamma_i^{\rm post} = \left(\frac{\partial \mathcal{G}}{\partial \x}^{\top} \Gamma_i^{-1} \frac{\partial \mathcal{G}}{\partial \x} + \Sigma^{-1}\right)^{-1}$. 
%For our specific problem, $\y_i \in \mathbb{R}^2$ and when the state space is high dimensional, $\frac{\partial \mathcal{G}}{\partial \x}^{\top} \Gamma_i^{-1} \frac{\partial \mathcal{G}}{\partial \x}$ might not be very informative. In such scenarios, one can approximate the posterior covariance with the prior covariance.
\begin{algorithm}
	\renewcommand{\algorithmicrequire}{\textbf{Input:}}
	\renewcommand{\algorithmicensure}{\textbf{Output:}}
	 \caption{Metropolis-Hastings algorithm to sample PDF $\pi$}
	\label{alg:MHalgorithm}
	\begin{algorithmic}[1]
		\Require Initial guess $x_1 \in \Omega$, $\pi$ (target distribution), and $Q$ (proposal distribution)
		\Ensure Samples from  $\pi$ namely,  $x_i$ for $i = 1,2, \dots$
                \Statex
                \State  Initialize $x_1$
                \For{$i = 1,2, \dots$}
                        \State Sample $\mathbf{z}$ from the proposal distribution, $Q( x_i, \mathbf{z})$
                        \State Evaluate $\displaystyle \alpha = \min \left(1,\frac{\pi(\mathbf{z}) Q(\mathbf{z}, x_i) }{\pi(x_i)  Q(x_i, \mathbf{z})} \right)$
                        \State Draw $ s \sim \mathcal{U}(0,1)$
                        \If{$\alpha > s$}
                                \State Accept: set $x_{i+1} = \mathbf{z}$
                        \Else 
                                \State Reject: set $x_{i+1} = x_{i}$
                        \EndIf
                \EndFor
	\end{algorithmic}
\end{algorithm}

\subsection{Markov chain Monte Carlo}
The Metropolis-Hastings (M-H) algorithm \cite{Metropolis_1953, Hastings_1970} is an MCMC method that employs a proposal density function ($Q$) at each sample point in $\Omega$ to generate a proposed sample point. This sample point is then rejected or accepted based on the M-H criterion ($\alpha$ in Algorithm \ref{alg:MHalgorithm}). The M-H algorithm is described in Algorithm \ref{alg:MHalgorithm} \cite[Section 3.6.2]{Kaipio_2006B}.
The performance of MCMC algorithms depends heavily on how close the proposal distribution is to the target distribution. A number of different MCMC algorithms exist, the distinguishing feature being the manner in which the sample points are proposed and accepted (or rejected). See, for example, \cite{Liu_2008B, Gilks_1995B, Gelfand_1990, Geman_1987}. In this paper, we use the delayed rejection adaptive Metropolis (DRAM) MCMC algorithm \cite{Haario_2006}.
\subsubsection{DRAM MCMC} \label{sec:DRAMMCMC}
     DRAM combines two ideas:  delayed rejection (DR) and adaptive Metropolis (AM) algorithms. Here, we describe DR, AM, and their combination.

     \subsubsection{Delayed rejection} DR is a strategy that is employed to improve the performance of the M-H algorithm. Unlike the M-H algorithm, which employs a single proposal density, DR uses a hierarchy of proposal densities. Suppose the current position of the Markov chain is $\Xi_n = \xi $ and that a candidate move $\Lambda_1$ is generated from the proposal distribution $q_1(\xi, \cdot)$. This proposal is accepted with probability 
     \begin{align*}
        \alpha_1(\xi, \lambda_1) = \min \left(1, \frac{\pi(\lambda_1) q_1(\lambda_1, \xi)}{\pi(\lambda_1) q_1(\xi, \lambda_1)}\right)\,.
     \end{align*}
     In the case of a rejection, the M-H algorithm retains the same position $\xi$. On the other hand, the DR algorithm instead  proposes a second move $\Lambda_2$. The second-stage proposal, $q_2(\xi, \lambda_1, \cdot)$, depends on the current position $\xi$ and  on the recently proposed and rejected move. The second-stage proposal is accepted with probability
     \begin{align}\label{eqn:acceptance_criterion}
        \alpha_2(\xi, \lambda_1, \lambda_2) = \min \left(1, \frac{\pi(\lambda_2)q_1(\lambda_2, \lambda_1) q_2(\lambda_2, \lambda_1, \xi)(1-\alpha_1(\lambda_2, \lambda_1))}{\pi(\xi)q_1(\xi, \lambda_1) q_2(\xi, \lambda_1, \lambda_2)(1-\alpha_1(\xi, \lambda_1))}\right)\,.
     \end{align}
     This process of delaying rejection can be iterated over a fixed number of stages. Alternatively, one can  use a biased coin to guide whether to move to a higher-stage proposal or not. We refer  interested readers to \cite{Haario_2006, Tierney_1999} for more details about the DR algorithm.

     \subsubsection{Adaptive Metropolis}
The AM MCMC algorithm constructs a proposal distribution adaptively by using the existing elements of the Markov chain. The basic idea is to use the sample path of the Markov chain to ``adapt'' the covariance matrix for a Gaussian proposal distribution. For example, after an initial period of nonadaptation, one can set the Gaussian proposal to be centered at the current position of the Markov chain, $\Xi_{n}$. That is, the covariance is set to $C_n  = s_d \textrm{Cov} (\Xi_0, \cdots, \Xi_{n-1}) + s_d \epsilon I_d\,$,
where $s_d$ is a parameter that depends only on the dimension of the state space on which the target probability distribution is defined. The quantity $\epsilon > 0$ is typically chosen to be a small constant, and $I_d$ is an identity matrix of appropriate dimensions. Before the start of the adaptation period, a strictly positive definite covariance $C_0$ is chosen according to a priori knowledge. Let index $n_0 > 0$ define the length of the nonadaption period. Then 
\begin{align}\label{eqn:AMCovariance}
        C_n =
        \begin{cases}
                C_0\,, \quad  & n  \leq n_0\,,\\
                s_d \textrm{Cov} (\Xi_0, \cdots, \Xi_{n-1}) + s_d \epsilon I_d\,, \quad & n > n_0 \,.
        \end{cases}
\end{align}
A recursive procedure allows us to update the covariance of the proposal distribution efficiently. For more details  see \cite{Haario_2001, Haario_2006}.
\subsubsection{Combining DR and AM} The success of the DR algorithm depends on a proposal in at least one of the stages being calibrated close to the target distribution. The AM algorithm attempts to calibrate the proposal distribution as the sample path of the Markov chain grows. The DRAM algorithm \cite{Haario_2006} combines these two strategies. The DRAM version deployed in this paper combines $m$ stages of DR with adaptation. The process can be summarized as follows:
\begin{itemize}
\item The proposal ($Q$ in Algorithm \ref{alg:MHalgorithm}) at the first of the $m$ stages is adapted  as described in equation \eqref{eqn:AMCovariance}. The covariance $C_n^1$ of the proposal distribution is computed by using the sample path of the Markov chain. 
\item The covariance $C_n^i$ of the proposal for stage $i$, $(i = 2,\cdots, m)$ is computed as a scaled version of the first-stage proposal $C_n^i = \gamma_i C_n^1$\,. 
\end{itemize}
Both $m$ and $\gamma_i$ can be freely chosen. For our purposes, we use a MATLAB implementation of DRAM that is available online \cite{Haario_2006}.
\section{Constructing the importance biasing distribution}\label{sec:ibd}
We explained in \S \ref{sec:mcmc-based-sampling} that solving Bayesian inverse problems is an effective method for constructing $X_i$'s (preimages to observations $\y_i$). The Bayesian inverse problem that we wish to solve is described in \ref{eqn:posterior2_approximate}. We also explained how to choose the covariance $\Gamma_i$ for the likelihood in question. In \S \ref{sec:BayesianInverseProblem}, we described two approaches,  Laplace approximation at MAP and DRAM MCMC, to solve the Bayesian inverse problem. One can use either of these approaches to draw samples from the unnormalized distribution $\pi_{\rm post}^i$ (DRAM) or an approximation of it (MAP). These samples are used to approximate the preimages $X_i$. In \S \ref{sec:mcmc-based-sampling} (equations \eqref{eqn:probability_conditional} and \eqref{eqn:approximate_omega_u}) we mentioned that $P_T(u)$ can be computed by approximating the set $\Omega(u)$ and using the corresponding probability measure $\mu$. A more practical means to estimate $P_T(u)$ is to use the preimages to construct an IBD and use the IBD to estimate the probability $P_T(u)$ using importance sampling. 

Using the Laplace approximation of the posterior, one can draw samples from the approximate posterior by sampling from the distribution in \eqref{eqn:laplace_approximation}, and these samples can be used to estimate $P_T(u)$ using IS. DRAM MCMC, on the other hand, yields a Markov chain, and we denote the elements of the Markov chain drawn from the unnormalized distribution $\pi_{\rm post}^i$ by
\begin{align}
        \widehat{X}_i := \{\x_i^1, \x_i^2, \cdots \}\,.
\end{align}
Assuming there are $\ell$ samples in the chain, the elements of set $\widehat{X}_i$ can be thought of as samples from the Gaussian distribution with empirical mean 
\begin{align}
        \overline{\x}_i = \frac{1}{\ell} \sum_{k=1}^{\ell} \x_i^k
\end{align}
and empirical covariance
\begin{align}
        \overline{X}_i = \frac{1}{\ell - 1} \sum_{k=1}^{\ell} (\x_i^k -\overline{\x}_i)(\x_i^k -\overline{\x}_i)^{\top}\,.
\end{align}
If we use $N$ observations (see the discussion around \eqref{eqn:approximate_omega_u}), then we can approximate the IBD as the following Gaussian mixture:
\begin{align}\label{eqn:ibd}
        p^{\rm IBD} :=& \sum_{i=1}^{N} \, w_i \mathcal{N}(\overline{\x}_i, \overline{X}_i)\,, \\
       \nonumber \sum_{i=1}^N w_i =& \,1\,.
\end{align}
One of the obvious ways to choose $w_i$ is to assign equal weights to each component of the mixture. This is  effective if $N$ is small, because the observations mostly correspond to high-density regions. If $N$ is large, however, the $(y_i,t_i)$ samples could potentially be from low-density regions, too. In such a scenario, it would be prudent to set 
\begin{align}
        w_i \propto y_i \varphi_t(u,y_i) \left |_{t=t_i} \right.\,.
\end{align}
\subsection{Estimating $P_T(u)$}\label{sec:pu}We now have all the pieces necessary to estimate $P_T(u)$. Following the discussion from \S \ref{sec:is}, the importance sampling estimate of $P_T(u)$ can be written as
\begin{align}\label{eqn:is_estimate}
        P_T^{\rm IS}(u)(\widehat{\mathbf{x}}_0^1, \ldots, \widehat{\mathbf{x}}_0^M) = \frac{1}{M} \sum_{i=1}^M\,\mathbb{I}(\widehat{\mathbf{x}}_0^i)\psi(\widehat{\mathbf{x}}_0^i)\,,
\end{align}
where $\widehat{\mathbf{x}}_0^1, \ldots, \widehat{\mathbf{x}}_0^M$ are sampled from the biasing distribution $p^{\rm IBD}$ and $\mathbb{I}(\widehat{\mathbf{x}}_0^i)$ represents the indicator function given by
\begin{align}
        \mathbb{I}(\widehat{\mathbf{x}}_0^i) =
\begin{cases}
1\,, \displaystyle \qquad \sup_{0 \leq t \leq T} \mathbf{c}^{\top}\mathbf{x}(t, \widehat{\mathbf{x}}_0^i) \geq u\,, ~~t\in \lbrack 0, T \rbrack\,,\\
0\,, \displaystyle \qquad \sup_{0 \leq t \leq T} \mathbf{c}^{\top}\mathbf{x}(t, \widehat{\mathbf{x}}_0^i) < u\,, ~~t\in \lbrack 0, T \rbrack\,.
\end{cases}
\end{align}
Also, $\psi(\widehat{\mathbf{x}}_0^i)$ represents the importance weights. The importance weight for an arbitrary $\widehat{\mathbf{x}}_0^i$ is given by
\begin{align}
     \displaystyle   \psi(\widehat{\mathbf{x}}_0^i) = \frac{p(\widehat{\x}_0^i)}{p^{\rm IBD}(\widehat{\x}_0^i)}\,.
\end{align}
 The overall procedure to compute an estimate of $P_T(u)$ is summarized in \Cref{alg:MCMC_Based_Is}.
\begin{algorithm}
        	\renewcommand{\algorithmicrequire}{\textbf{Input:}}
        	\renewcommand{\algorithmicensure}{\textbf{Output:}}
        	 \caption{Algorithm to estimate $P_T(u)$}
        	\label{alg:MCMC_Based_Is}
        	\begin{algorithmic}[1]
        		\Require Dynamics \eqref{eqn:dynamics}, initial distribution of the state $p$, and excursion level $u$
        		\Ensure An estimate of $P_T(u)$
                        \Statex
                        \For{$i = 1,2, \dots$}
                                \State  Sample from $y \varphi_t(u,y)$ using DRAM MCMC algorithm described in \S \ref{sec:DRAMMCMC}to construct $\y_{i}$. Use the details given in \S \ref{sec:evaluatevarphi} to evaluate $y \varphi_t(u,y)$.
                                \State  Construct the likelihood by using the formula in equation \eqref{eqn:likelihood}. The covariance information can be constructed by using the approach in \S \ref{sec:covarianceGamma}.
                                \State Construct the posterior distribution by using the formula in equation \eqref{eqn:posterior2}.
                                \State Generate samples from approximate $\pi_{\rm post}^i$ by using either the Laplace approximation at MAP or the DRAM MCMC algorithm (details in \S \ref{sec:BayesianInverseProblem})
                                \State Use the samples obtained in the previous step to construct the IBD (details  in \S \ref{sec:ibd}, specifically equation \eqref{eqn:ibd}).
                        \EndFor
                        \Statex Use the formulae in \eqref{eqn:is_estimate} to obtain $P^{\rm IS}_T(u)$
        	\end{algorithmic}
\end{algorithm}
We call the approach that uses the MAP point and Laplace approximation around the MAP point to construct the IBD as {\textit{MAP-based IS}} and the approach that uses the MCMC chains to construct the IBD as {\textit{MCMC-based IS}}.
\subsection{Connection to approach based on LDT} We  mentioned earlier that LDT has been used in \cite{Dematteis_2018, Dematteis_2019} to estimate rare-event probabilities using large deviations as a tool. For a detailed treatment of large deviation theory, we refer the interested readers to \cite{Touchette_2011}. Loosely speaking, one can use large deviations to estimate $P_T(u)$ when $P_T(u) \rightarrow 0$ as $u \rightarrow \infty$. According to LDT,
\begin{align}
P_T(u) \asymp \exp(-  I(\x))\,,
\end{align}
where $\asymp$ indicates that the ratio of the logarithm 's right-hand side and the logarithm's left-hand side tends to one asymptotically, where 
\begin{align}
        I(\x) := \displaystyle \frac{1}{2}\, \underset{\x \in \Omega(u)}{\textrm{min}} \|\x - \overline{\x}\|^2_{\Sigma^{-1}}\,.
\end{align}

Intuitively this approach, introduced by Dematteis et al.~in \cite{Dematteis_2018, Dematteis_2019}, approximates the rare-event probability by determining the dominating point in $\Omega(u)$, and the relative precision of this estimate improves as $u$ increases. On the other hand, for given $u$, which is the case  we discuss here, even determining an error estimate for the large deviation approach is problematic in practice. While inspired by large deviation ideas \cite{adler2009random}, our approach goes further by approximating the distribution around the dominating point the distribution in an importance sampling approach that  produces an unbiased estimate of the sought-after probability. The empirical variance of the importance sampling approach gives an estimate of the error we make in our approach, something that is not accessible in a classical large deviation approach. 

\section{Numerical results}\label{sec:num_exp}
We demonstrate the application of procedure described in \S \ref{sec:mcmc-based-sampling} and \S \ref{sec:BayesianInverseProblem} for nonlinear dynamical systems excited by a Gaussian distribution. We use the Lotka-Volterra equations and the Lorenz-96 system as test problems.
\subsection{Lotka-Volterra system}\label{sec:nonlinear_Gaussian}
\begin{figure}
                        \begin{center}
                          \subfigure[Contours of $y\varphi_t$ for $u = 17$ ]{\includegraphics[width=0.49\linewidth]{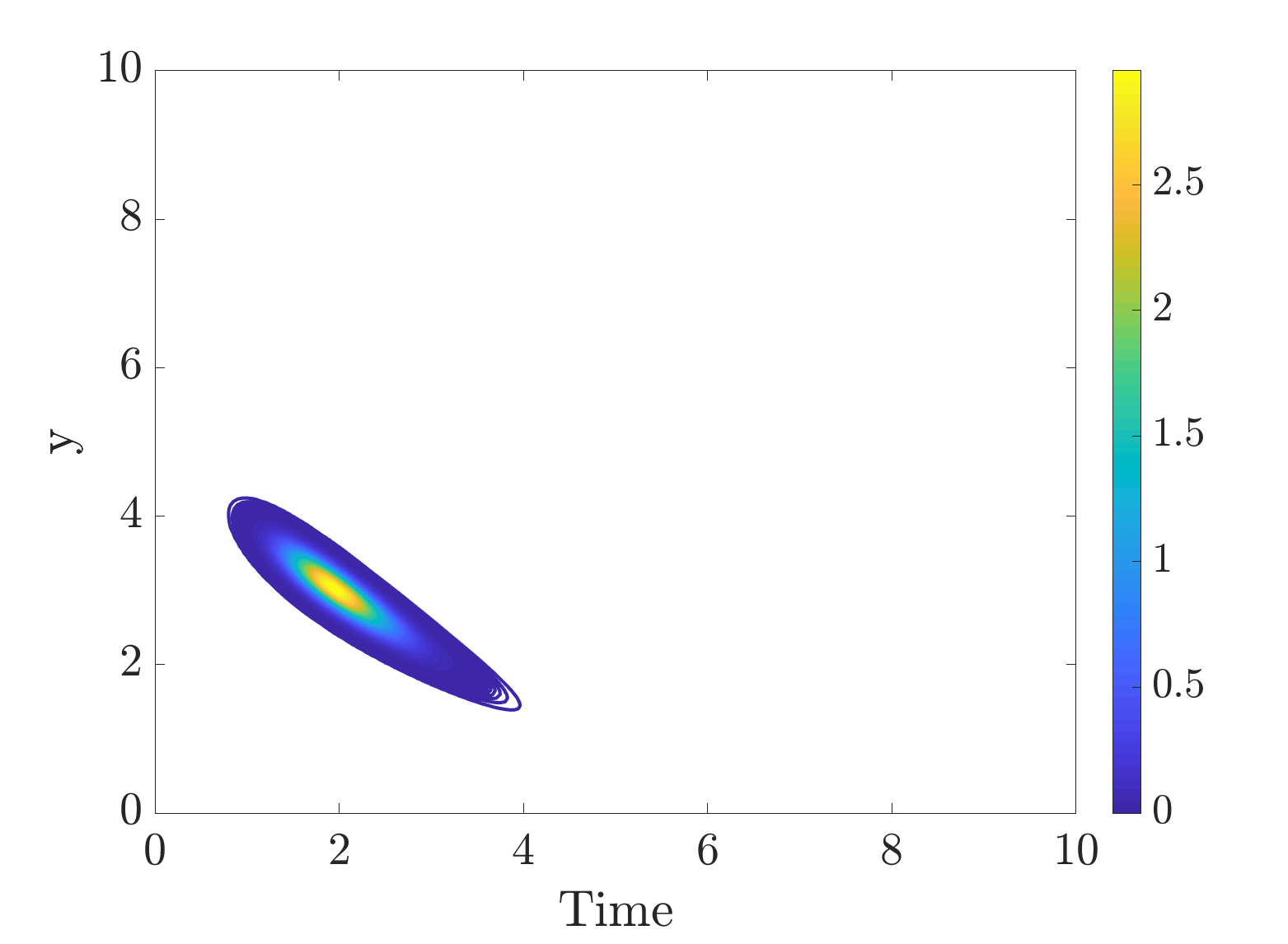}}
                          \subfigure[Samples drawn from $y\varphi_t$ using DRAM MCMC]{\includegraphics[width=0.49\linewidth]{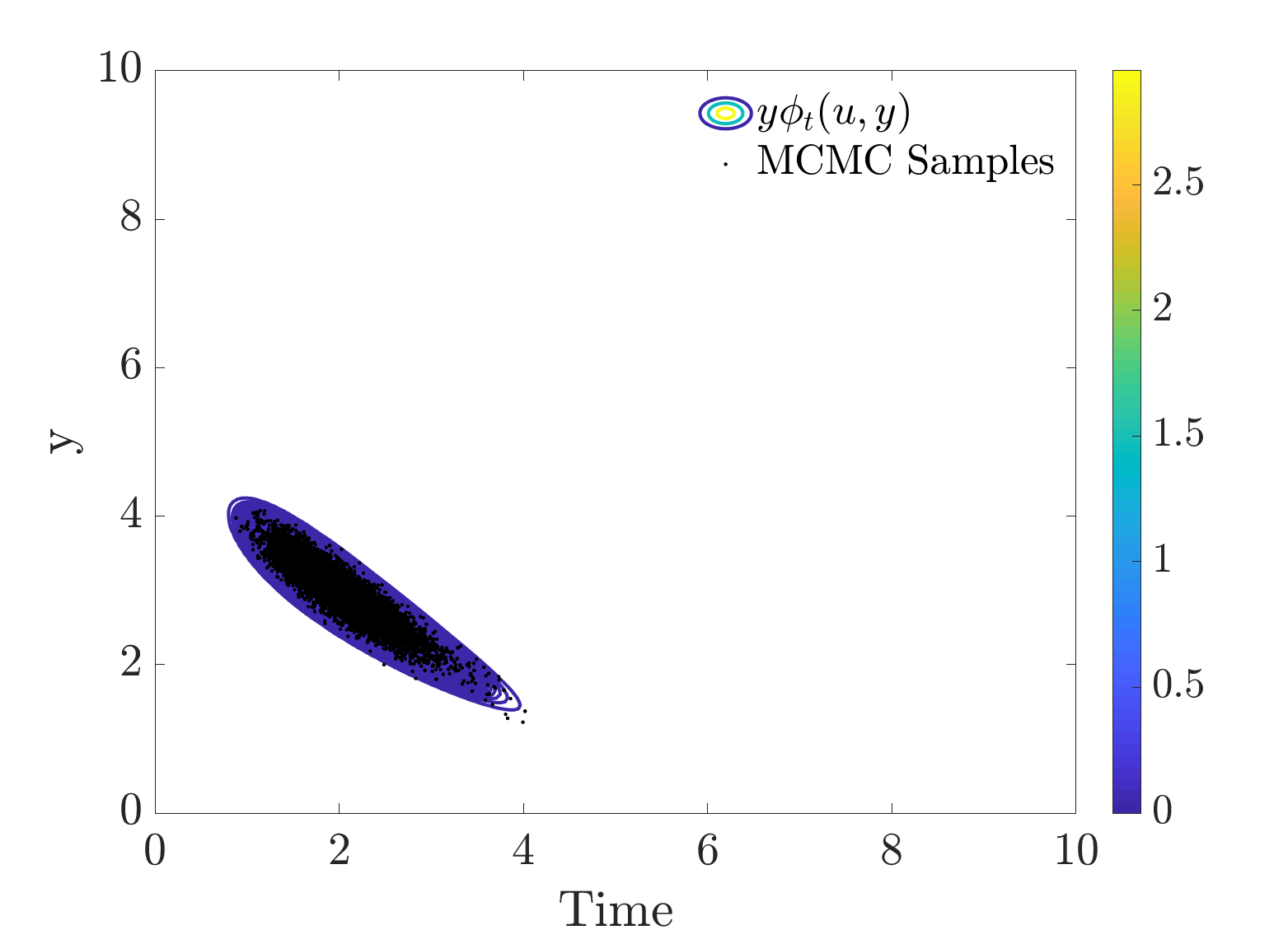}}
                          \end{center}
                         \caption{Left: Product of the derivative and the joint PDF of the state and its derivative for $u=17$. Right: Samples drawn from  $y\varphi_t$ using DRAM MCMC. These samples will be used to construct $\y_{\rm obs}^i$, which in turn will be used to construct $\pi_{\rm post}^{i}$ }
                         \label{fig:contoursofphi_U_20}
                \end{figure}

\begin{figure}
                        \begin{center}
                          \subfigure[Autocorrelation between samples drawn from $y\varphi_t$ ]{\includegraphics[width=0.55\linewidth]{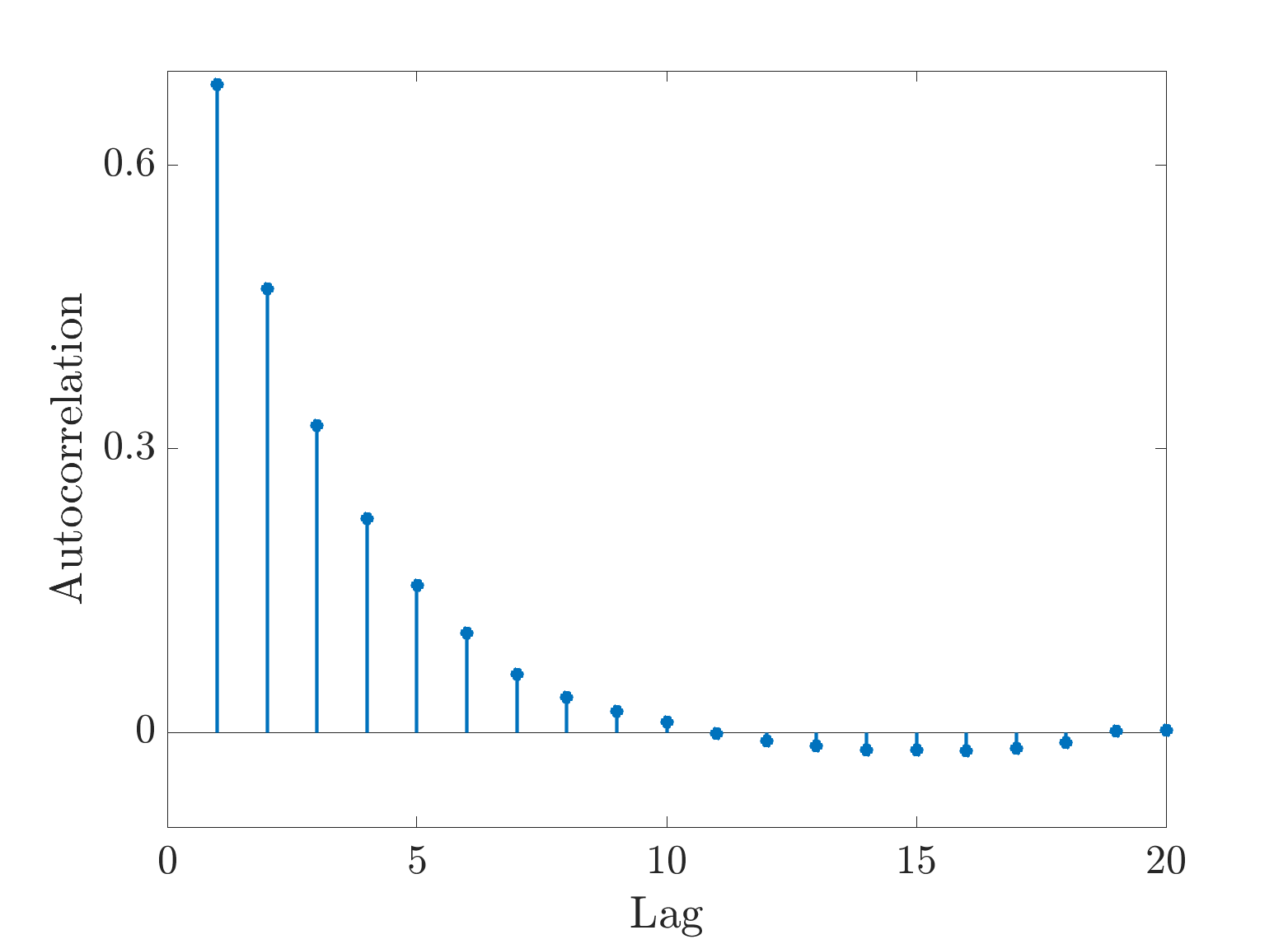}}
                         \caption{Autocorrelation function vs the lag for samples generated from $y\varphi_t$. We see that autocorrelation dies down to zero for a lag of $11$.}
                        \end{center}
                         \label{fig:Autocorrelation_Varphi}
                \end{figure}

\begin{figure}
                        \begin{center}
                          \subfigure[Samples from nominal and biasing distributions]{\includegraphics[width=0.85\linewidth]{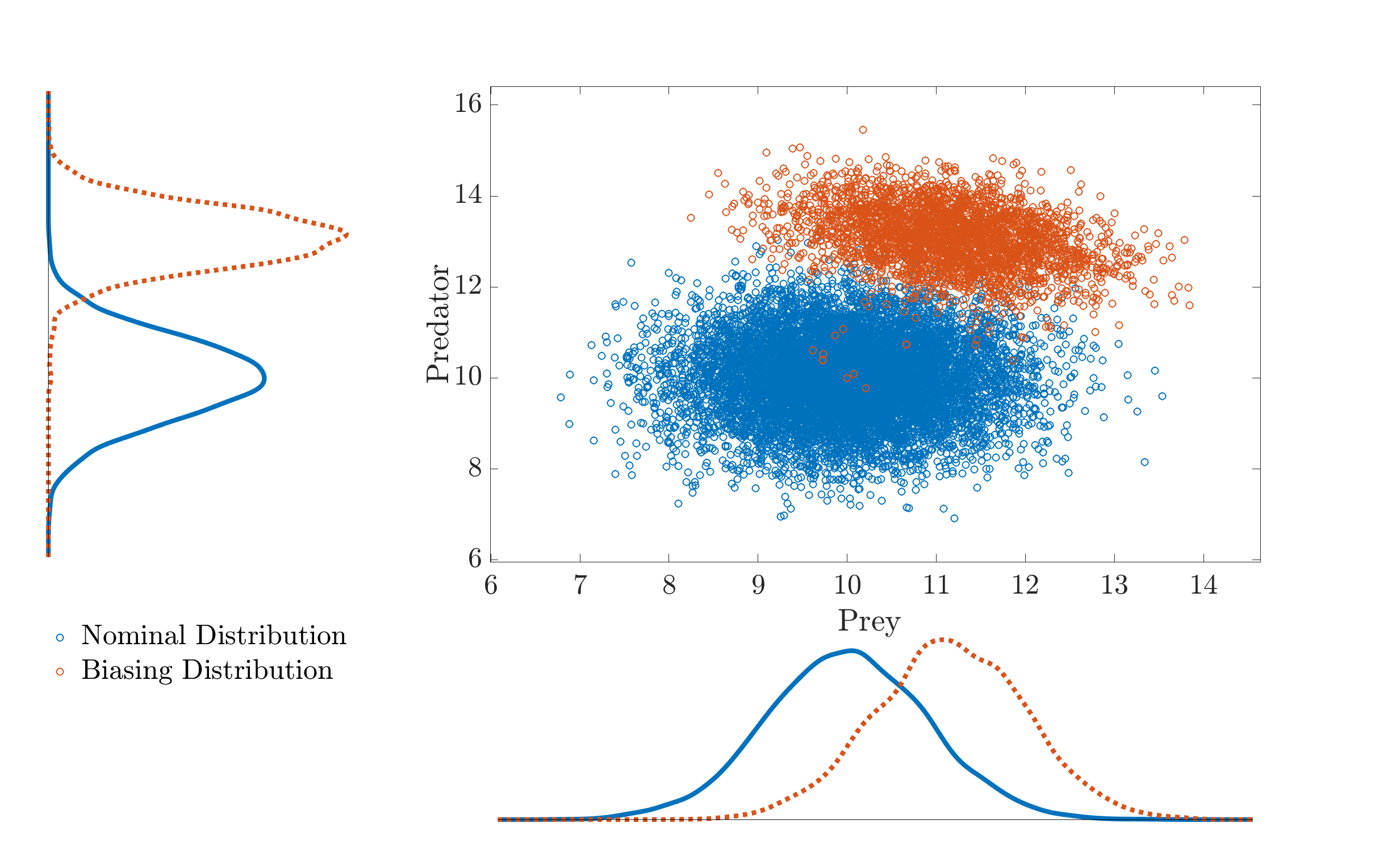}}
                         \caption{Samples from nominal and biasing distributions. The biasing distribution is constructed by using five observations.}
                        \end{center}
                         \label{fig:nominal_biasing}
                \end{figure}

\begin{figure}
                        \begin{center}
                          \subfigure[Convergence of MCMC-based IS and MCS for 2D Lotka-Volterra with Gaussian input]{\includegraphics[width=0.55\linewidth]{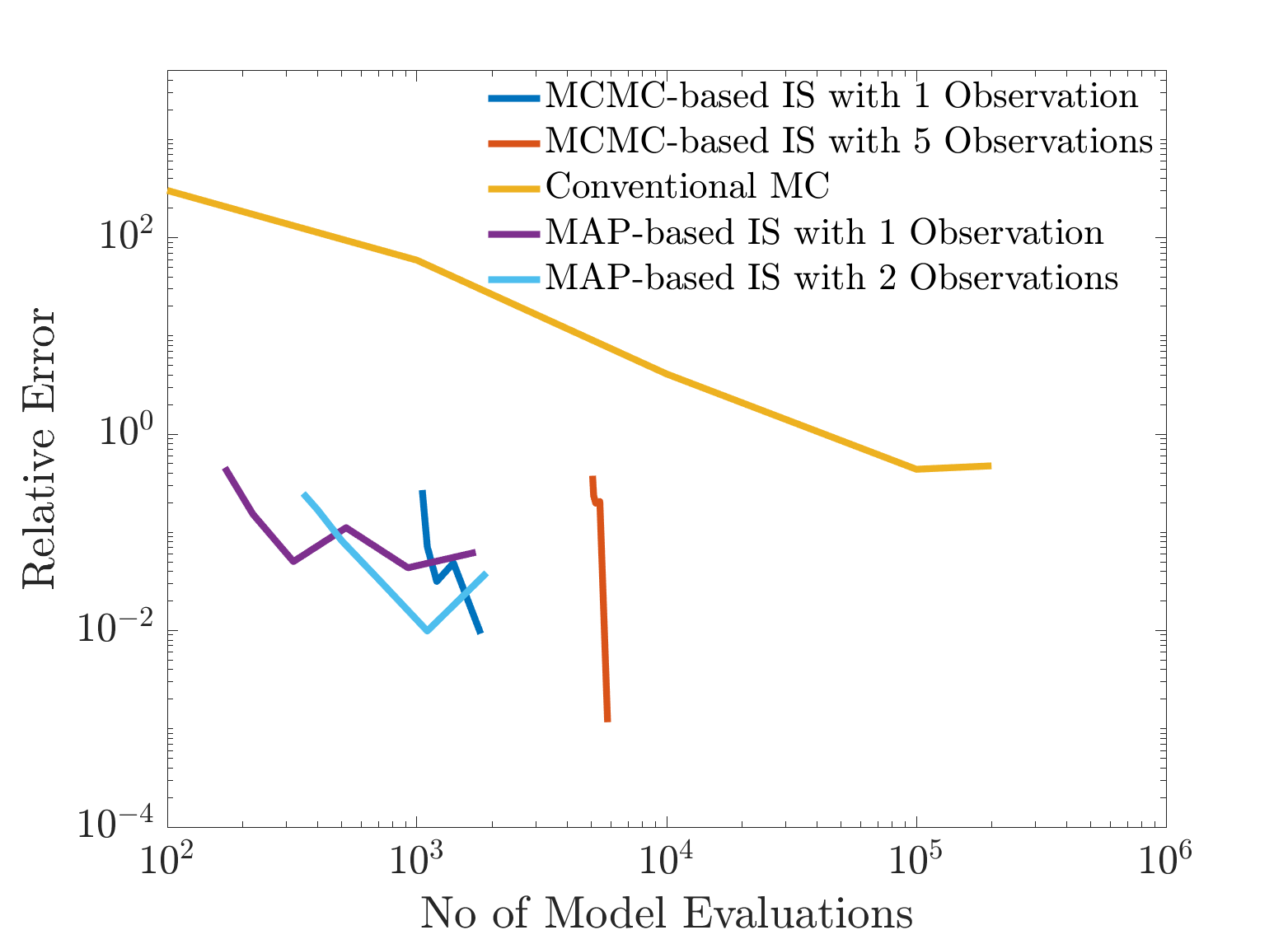}}
                         \caption{Convergence of MCMC-based IS, MAP-based IS, and MCS. The true proability here is $3.28\times 10^{-5}$. With about 1,000 model evaluations, MAP-based IS yields a 1\% accurate probability estimate. MCMC-based IS converges rapidly; with about 1,000 model evaluations, we see a fairly accurate estimate; and with about 5,000 samples, the accuracy of the estimate is much better.}
                        \end{center}
                         \label{fig:nonlinear_convergence}
                \end{figure}

The Lotka-Volterra equations, which are also known as the predator-prey equations, are a pair of first-order nonlinear differential equations and are used to describe the dynamics of biological systems in which two species interact, one as  predator and the other as  prey. The populations change through time according to the following pair of equations,
\begin{align}\label{eqn:nonlinear_example}
\displaystyle \frac{dx_1}{dt} = \alpha x_1 - \beta x_1x_2 \,, \\
\nonumber \displaystyle \frac{dx_2}{dt} = \delta x_1x_2 - \gamma x_2\,,
\end{align}
where $x_1$ is the number of prey, $x_2$ is the number of predators, and $\displaystyle \frac{dx_1}{dt}$ and $\displaystyle \frac{dx_2}{dt}$ represent the instantaneous growth rates of the two populations. We assume that the initial state of the system at time $t=0$ is a random variable that is normally distributed:
$$\x(0) \sim \mathcal{N}\left(\begin{bmatrix}10 \\ 10 \end{bmatrix}, 0.8\times I_2\right)$$. We are interested in estimating the probability of the event $P(\mathbf{c}^{\top}\x \geq u)$, where $\mathbf{c} = \begin{bmatrix} 0 \\ 1 \end{bmatrix}$, $t \in [0,10]$, and $u = 17$. The first step of our solution procedure involves sampling from $y \varphi_t(u,y)$ to generate observations $\y_i$. We  linearize the dynamical system about the mean of the distribution of $\x_0$ equation  \eqref{eqn:linearizedsystem} and express $\varphi_t(u,y)$ as a function of $t$ and $y$ as described by equation \eqref{eqn:gpvarphi}. We can compute $y \varphi_t(u,y)$  as shown in equation \eqref{eqn:compute_varphi}. We use the DRAM MCMC method to generate samples from $y \varphi_t(u,y)$; to minimize the effect of the initial guess on the posterior inference, we use a burn-in of 1,000 samples. \cref{fig:contoursofphi_U_20} shows the contours of $y \varphi_t(u,y)$ and samples drawn from it by using DRAM MCMC. \cref{fig:Autocorrelation_Varphi} shows the autocorrelation between the samples drawn by using DRAM MCMC from $y \varphi_t(u,y)$, and we see that the autocorrelation dies down to zero for a lag of $11$; choosing every eleventh sample gives us independent samples that are in turn used to form $\y_i$. The next step in our solution procedure is to construct $\pi_{\rm post}^{i}$ that approximate the preimages of $\y_i$. We use the procedure described in \S \ref{sec:mcmc-based-sampling} to form an unnormalized posterior distribution that uses $\y_i$. Subsequently, either MAP-based IS or MCMC-based IS can be used to estimate $P_T(u)$. For the MAP-based IS, we first solve the optimization problem in \eqref{eqn:ip_laplace} using a gradient-based optimization algorithm (for example, LBFGS). We use the inverse of the Hessian at the MAP point to approximate the covariance of the posterior as in equation \eqref{eqn:laplace_approximation}. 
For the MCMC-based IS, we use DRAM MCMC as described in \S \ref{sec:DRAMMCMC} to sample from the posterior distribution. To minimize the effect of initial guess on the posterior samples, we use a burn-in of 500 samples. We use these samples to construct the IBD as described in \S \ref{sec:ibd}. We test our algorithm by constructing $p^{\rm IBD}$ using different numbers of observations. \cref{fig:nominal_biasing} shows samples drawn from $p$, $p^{\rm IBD}$, and the corresponding marginal densities. We see that the samples generated from the IBD are predominantly from the tails of $p$.
\cref{fig:nonlinear_convergence} compares the relative accuracies of conventional MCS and MCMC-based IS algorithms. We use the Monte Carlo estimate obtained using 10 million samples as a proxy for the truth. We test the accuracy of the IBD constructed with $1$ and $5$ observations (see the discussion around \eqref{eqn:approximate_omega_u} for the definition of the number of observations). Constructing an IBD with $5$ observations involves more work because the MCMC DRAM has to be run with $5$ different unnormalized posterior distributions, involving about 5,000 model evaluations just to construct $p^{\rm IBD}$ and a further 800 model runs to estimate $P_T(u)$. We note that executing the MCMC DRAM with 5 different observations completely independent of one another  can be run in parallel. On the other hand, constructing $p^{\rm IBD}$ with a single observation requires 1,000 model runs and a further 800 model runs to estimate  $P_T(u)$. As \cref{fig:nonlinear_convergence} indicates, we get a more accurate estimate (one order of magnitude) for extra work performed with $5$ observations. For most practical purposes, however, an accuracy of 1\% that is obtained with $p^{\rm IBD}$ constructed from a single observation is sufficient. 

\begin{figure}
                        \begin{center}
                          \subfigure[Confidence intervals for MCMC-based IS with 1 observation]{\includegraphics[width=0.49\linewidth]{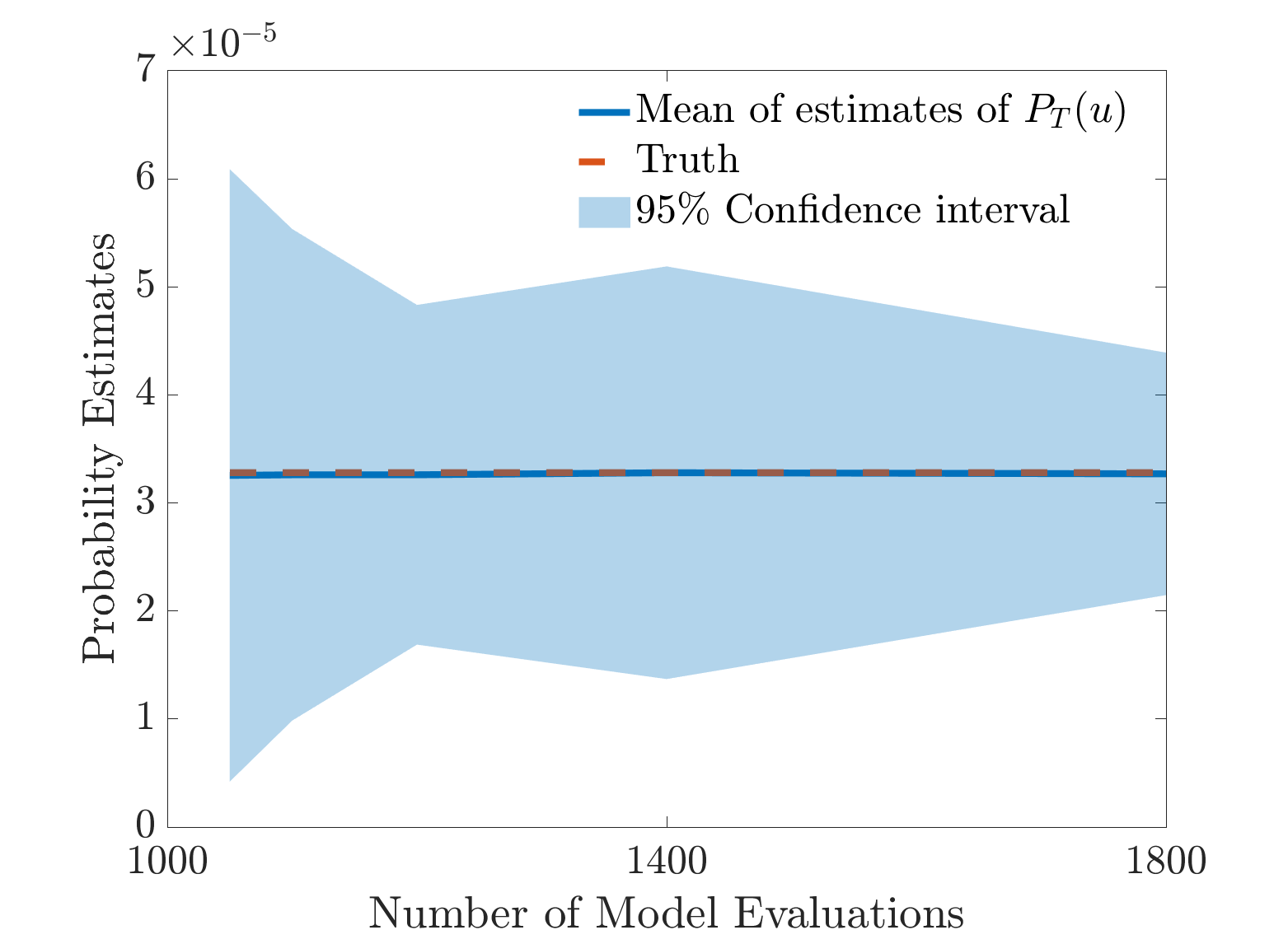}}
                          \subfigure[Confidence intervals for MCMC-based IS with 5 observations]{\includegraphics[width=0.49\linewidth]{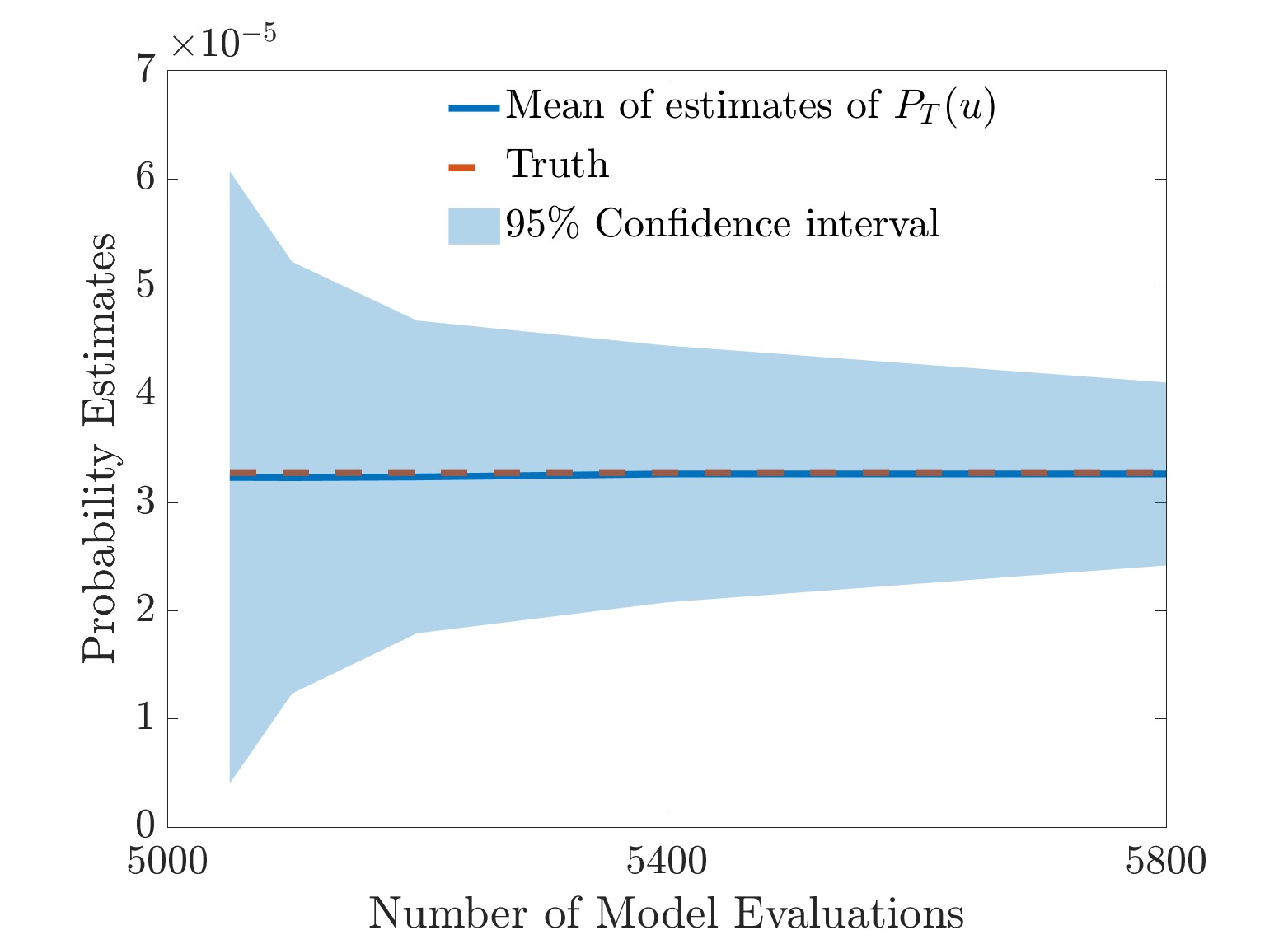}}
                          \subfigure[Confidence intervals for MAP-based IS with 1 observation]{\includegraphics[width=0.49\linewidth]{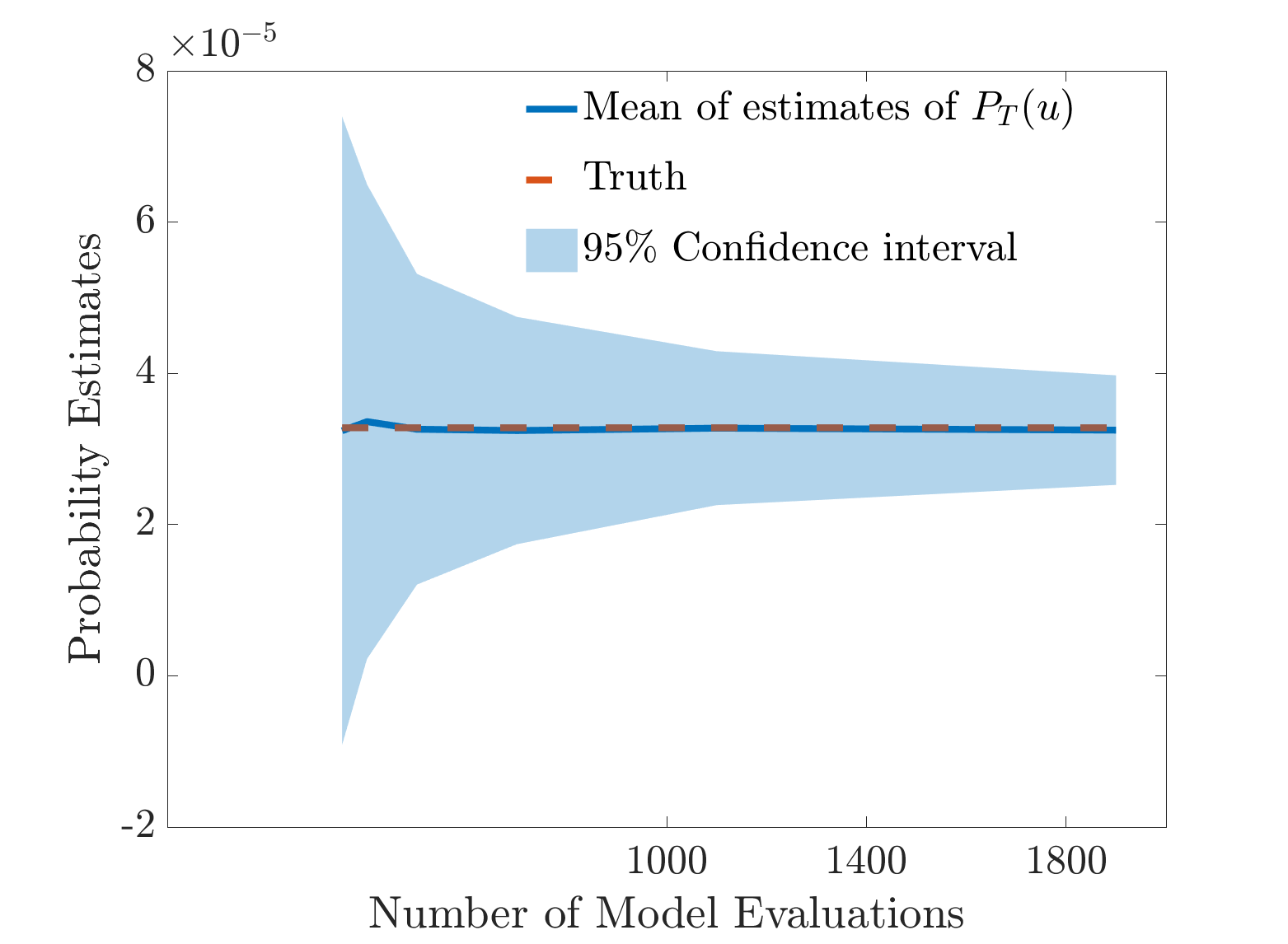}}
                         \caption{The 95\% confidence intervals of MAP-based IS with one observations and MCMC-based IS with one and five observations. The confidence interval for MAP-based IS with 1,800 model evaluations is sharper than that of MCMC-based IS with more than 5,000 model evaluations.}
                        \end{center}
                         \label{fig:confidence_intervals}
                \end{figure}

A certain amount randomness exists in almost every step of our algorithm. For a fair comparison, instead of reporting just a single error plot, we also report the confidence intervals of the estimates. We execute the algorithm \ref{alg:MCMC_Based_Is} 1,000 times and estimate the 95\% confidence intervals based on the results obtained with these runs. \cref{fig:confidence_intervals} shows the mean of the estimate, the truth, and the 95\% confidence intervals for $p^{\rm IBD}$ constructed with one and five observations. When $p^{\rm IBD}$ is constructed with one observation, our algorithm yields an estimate that is within 32\% of the actual value of excursion probability with 95\% probability. We note that to obtain an estimate within 32\% error, MCS requires $\mathcal{O}(10^5)$ model evaluations, whereas our method requires $\mathcal{O}(10^3)$ model evaluations. The estimates obtained with $p^{\rm IBD}$ constructed with five observations are sharper; that is, the confidence intervals are narrower. With $\mathcal{O}(5 \cdot 10^3)$ model evaluations, our algorithm yields an estimate that is within 25\% error with 95\% probability. To obtain the same level of accuracy, MCS will require $\mathcal{O}(5 \cdot 10^5)$ model evaluations. 

                % \begin{figure}
                %         \begin{center}
                %           \subfigure[Convergence of MCMC-based IS and MCS]{\includegraphics[width=0.6\linewidth]{Figures_Paper/CombinedErrorsWithDifferentH.png}}
                %          \caption{Convergence of MCMC-based IS and MCS for different choices of $h$. The true proability here is $3.28\times 10^{-5}$. Although the probability estimates are reasonably accurate for all choices of $h$, the convergence is most interpretable for $h$ = 1.}
                %         \end{center}
                %          \label{fig:Bhat_Hs}
                % \end{figure}
% \subsubsection{Number of observations to approximate $\widehat{\mathcal{B}}_u$}As discussed in section \ref{sec:computational_cons} the computational cost for constructing the IBD depends on the number of observations used to approximate $\widehat{\mathcal{B}}_u$.  Figure the accuracy of the probability estimates using MCMC based importance sampling using 100 samples with different approximations of $\widehat{\mathcal{B}}_u$ (based on the number of observations used). 
% \begin{figure}
%         \begin{center}
%           \subfigure[Accuracy of MCMC-based IS]{\includegraphics[width=0.45\linewidth]{Figures_Paper/Error_Vs_Obs.png}}
%          \caption{Accuracy of the estimates of $P_u$ for different approximations of $\widehat{\mathcal{B}}_u$ using 100 samples for importance sampling. ``True" probability is $3.28 \times {10}^{-5}$. Even with small number of observations, the estimates are accurate.}
%         \end{center}
%          \label{fig:Error_Vs_Obs}
% \end{figure}
%
% Commenting out Reusing the IBD part.
% \input{Reusing_The_IBD.tex}
%
\subsection{Lorenz-96 system}
\begin{figure}
                \begin{center}
                \subfigure[MAP-based IS with Lorenz-96 for Lorenz96 in 100 dimensions with Gaussian input]{\includegraphics[width=0.55\linewidth]{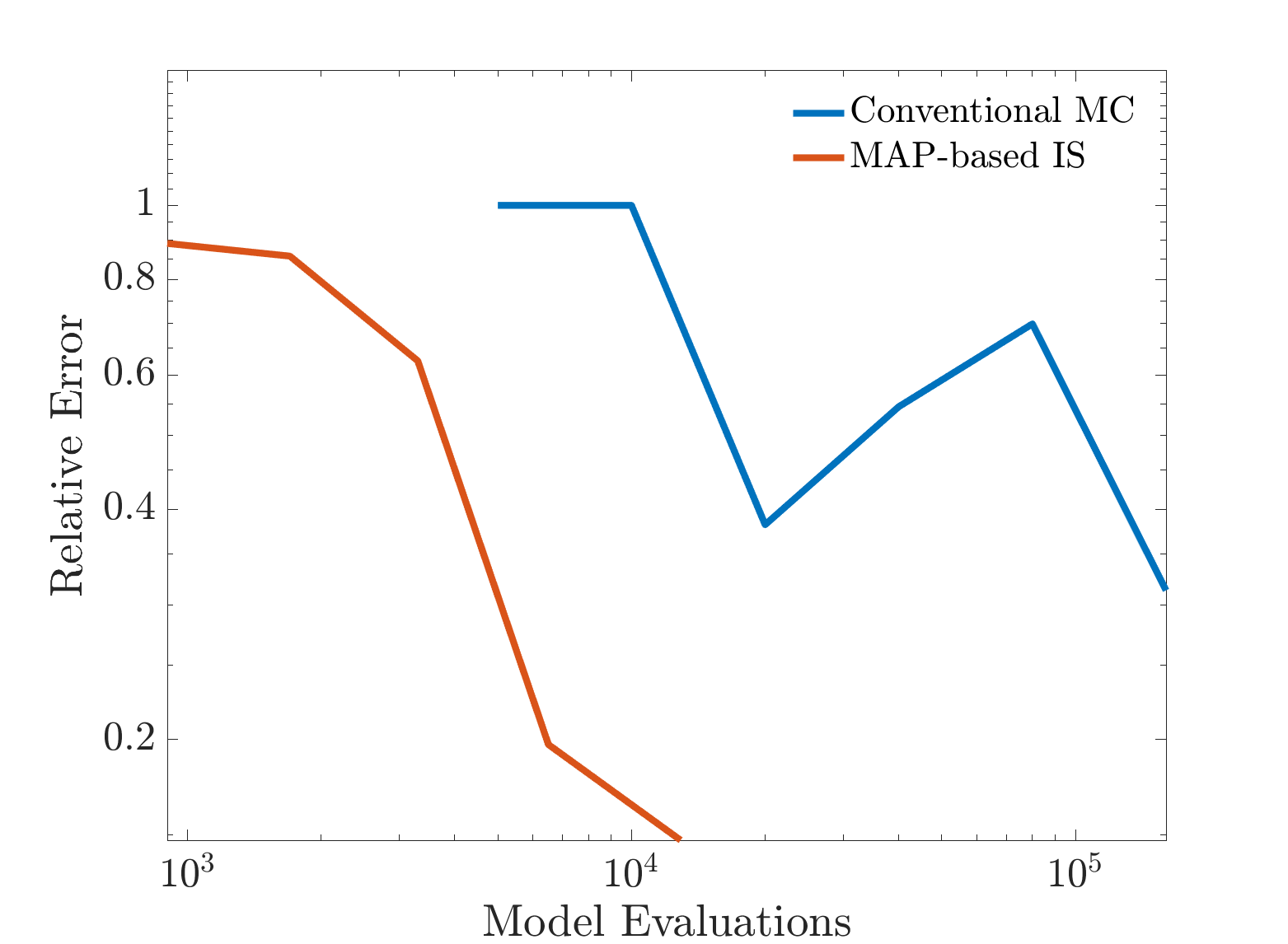}}
                  \end{center}
                  \vspace{0.05in}
                \caption{Convergence of MAP-based IS and MCS. The true proability here is $8.09\times 10^{-5}$. MAP-based IS converges rapidly. With about 5e3 model evaluations, we see a fairly accurate estimate; and with about 1e4 samples, the accuracy of the estimate is much better.}
                 \label{fig:estimate_100d}
        \end{figure}

Lorenz-96 system \cite{Lorenz_1996} is a one dimensional atmospheric model used to study the predictability of the atmosphere and weather forecasting. The system can be interpreted as atmospheric waves traveling around a circle of constant latitude. The equation of the dynamical system is 
\begin{align}
        \displaystyle \frac{\mathrm{d}x_i}{\mathrm{d}t} = x_{i-1}(x_{i+1}-x_{i-2}) - x_i + F\,, \quad i=1, \dots, n > 3\,,
\end{align}
with the periodic boundary conditions $x_{i+n} = x_i$. The Lorenz-96 model has been used as a test problem for data assimilation algorithms, subgrid scale parameterizations, and predictability of extreme waves \cite{Rao_2016, Crommelin_2008, Sterk_2017}. We demonstrate our methodology for a 100-dimensional Lorenz-96 system with $F = 3$. We are interested in estimating the probability of the event $P(\mathbf{c}^{\top}\x \geq u)$, where $\mathbf{c} = \begin{bmatrix} 1 \\ 0 \\ \vdots \\ 0 \end{bmatrix}$, $t \in [0,2]$, and $u = 6$. We use the Monte Carlo estimate obtained with 10 million samples as a proxy for the truth. For the Lotka-Volterra system, we observe that MAP-based IS gives  probability estimates (and confidence intervals) that are comparably accurate to those of MCMC-based IS with much fewer model evaluations. Hence, for the Lorenz-96 system, we demonstrate the results with MAP-based IS only. \Cref{fig:estimate_100d} demonstrates the convergence of the MAP-based IS and conventional MCS approaches. We observe that MAP-based IS achieves the same level of accuracy as does MCS with about 1\% of the computational effort.

        \begin{figure}
                \begin{center}
                  \subfigure[MCMC-based IS and MAP-based IS for the Lotka-Volterra system (2D) with an uniform excitation]{\includegraphics[width=0.55\linewidth]{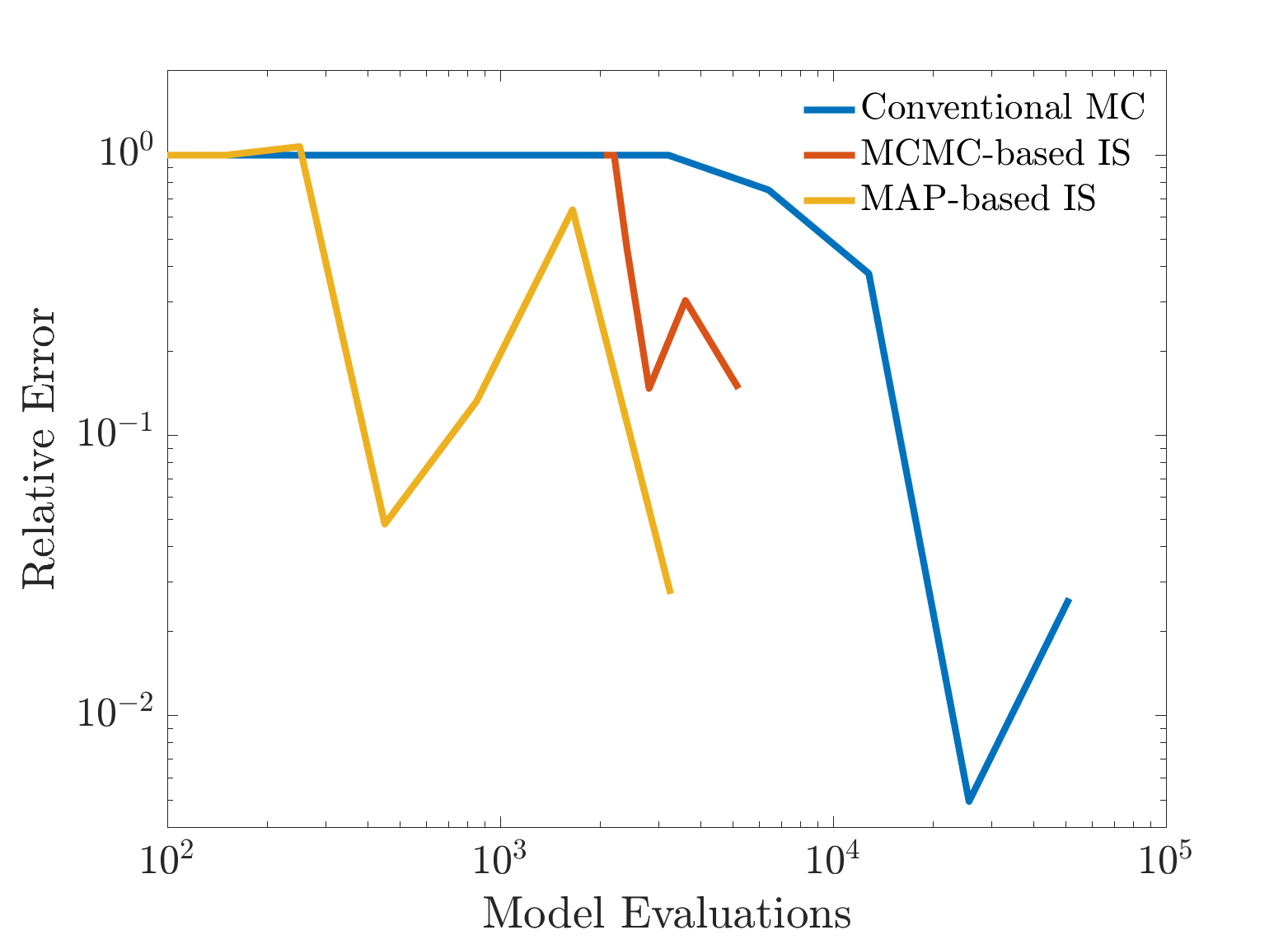}}
                  \caption{Convergence of MAP-based IS, MCMC-based IS, and MCS for the Lotka-Volterra system with a uniform excitation. The true probability here is $6.281\times 10^{-4}$. The convergence is not as smooth as it is for a Gaussian excitation, and we attribute the cause to the edge effects of a uniform distribution.}
                \end{center}
                \label{fig:non_Gaussian_LV}
        \end{figure}
        \begin{figure}
                \begin{center}
                \subfigure[Convergence of MAP-based IS for the Lorenz96 (100D) system with a uniform excitation]{\includegraphics[width=0.55\linewidth]{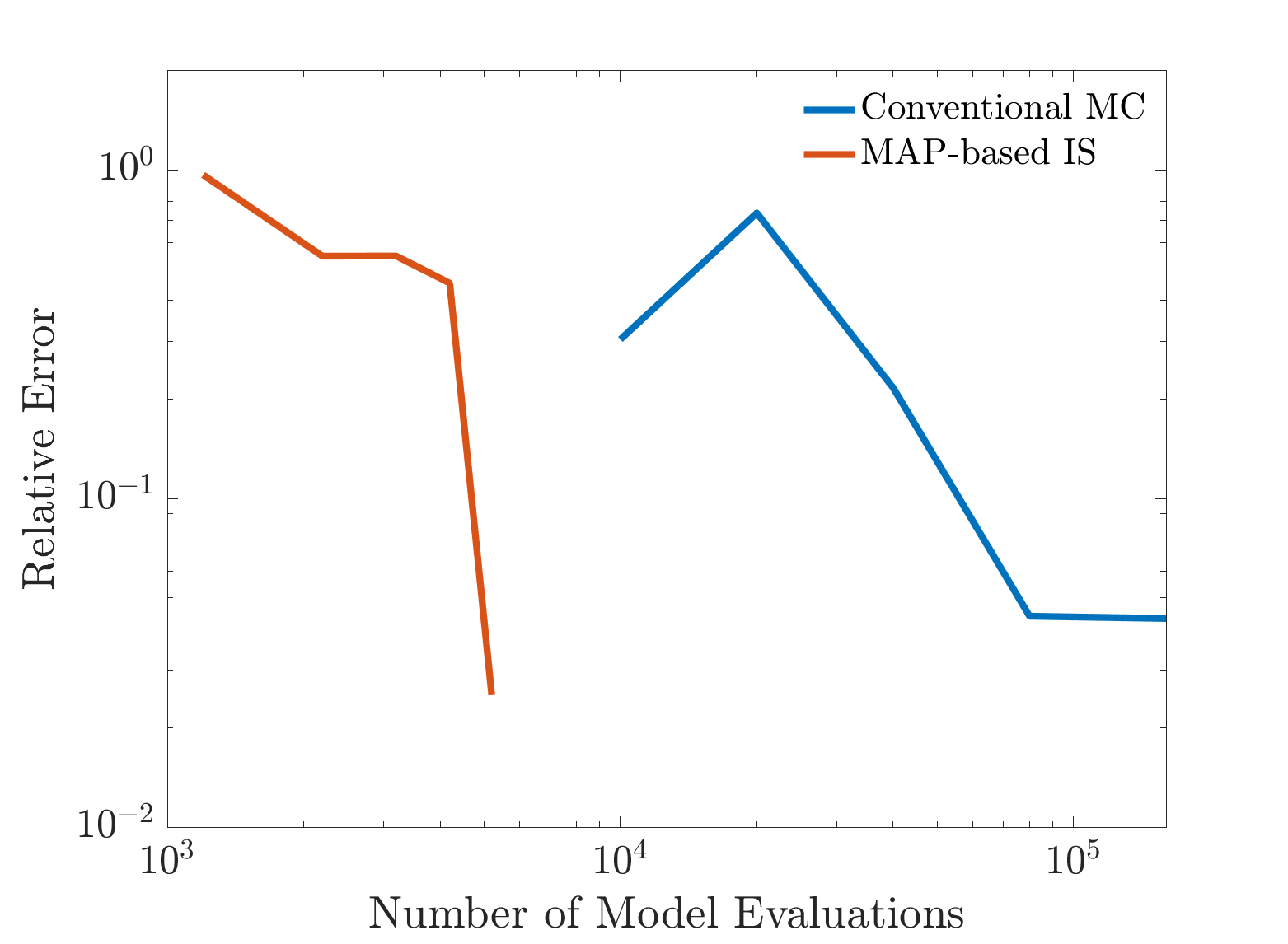}}
                \caption{Convergence of MAP-based IS and MCS for the Lorenz96 with a uniform excitation. The true proability here is $1.438\times 10^{-4}$. The probability estimates are as accurate as the MCS with less than 1\% of the computational cost.}
               \end{center}
                \label{fig:non_Gaussian_Lorenz}
        \end{figure}
                \subsection{Experiments with non-Gaussian excitation} All the numerical experiments discussed until now have been excited by a Gaussian random variable. In many scenarios, however, the dynamical system might be excited by a non-Gaussian random variable. In such a scenario, the method introduced in this paper cannot be used directly. As discussed  in \S \ref{sec:approxB}, for non-Gaussian processes $y\varphi_t(u,y)$ cannot be evaluated analytically. Hence, to overcome this limitation, we used the method of moments to approximate the non-Gaussian excitation. Subsequent steps of the computational procedure remain unaltered. We demonstrate this method by using a uniform distribution to excite the dynamical systems. For the Lotka-Volterra system described in \S \ref{sec:nonlinear_Gaussian}, we use $\x(0) \sim \begin{bmatrix} 9 \\ 9 \end{bmatrix} + 3.6 \times   \begin{bmatrix} \mathcal{U}(0,1) \\ \mathcal{U}(0,1) \end{bmatrix}$. For the Lorenz system we use $\x(0) \sim 0.5 + 3 \times \mathcal{U}[0,1]^{100}$. We approximate $\x(0)$ by using the method of moments and use this approximation in subsequent steps of the computational procedure. \Cref{fig:non_Gaussian_LV} demonstrates the results with both MAP-based IS and MCMC-based IS for the Lotka-Volterra system. We use the MC estimate with 10 million samples as a proxy for the truth. We observe that the convergence is not smooth, which could be potentially due to the edge effects and lack of a ``tail" in the uniform distributions. \Cref{fig:non_Gaussian_Lorenz} demonstrates the results with MCMC-based IS for the Lotka-Volterra system. We use the MC estimate with 10 million samples as a proxy for the truth. We observe that the probability estimates are accurate in this case and even the convergence is smooth. We obtain a probability estimate that is as accurate as the MCS estimate with about 1\% of the computational cost.
\section{Conclusions and future directions}\label{sec:conc}
In this paper, we have developed a novel method that uses excursion probability theory in conjunction with MCMC to estimate probabilities of rare and extreme events. The method uses Rice's formula to construct an IBD by means of Bayesian inference. As we have demonstrated in \S \ref{sec:num_exp}, the method gives accurate estimates of the probability with a small number of evaluations of the associated computational model. The results show that the algorithm obtains an estimate that is as accurate as conventional MCS with about 1\% of the computational effort. We note that the most expensive part of our algorithm is constructing the IBD. The computational burden can be alleviated by carrying out the construction of MCMC chains in parallel when multiple observations are used for constructing the IBD. For the MAP-based IS, we need to solve an optimization problem (or solve the same number of optimization problems as the number of observations used for constructing the IBD).

Currently, the method is feasible for problems with parameter dimension of $\mathcal{O}(100)$. Scaling this method to dimensions of $\mathcal{O}(1000)$ can be challenging; we will explore using surrogate models to alleviate the computational burden in  constructing the IBD. One could also use surrogate models of different fidelities and adaptively choose an appropriate model for this step. This approach could result in significant computational savings. The surrogate models of different fidelities could also be used in a multilevel framework to construct the Markov chain.

\section*{Acknowledgements}
We thank Emil Constantinescu and Charlotte Haley, our colleagues at Argonne for useful discussions during various stages of this work. Vishwas Rao also thanks Nick Alger from University of Texas, Austin for helpful comments on the paper.
\bibliography{ExcursionProbability}
%\vfill
% \begin{flushright}
% \scriptsize \framebox{\parbox{3.2in}{Government License
% The submitted manuscript has been created by UChicago Argonne, LLC,
% Operator of Argonne National Laboratory (``Argonne").  Argonne, a
% U.S. Department of Energy Office of Science laboratory, is operated
% under Contract No. DE-AC02-06CH11357.  The U.S. Government retains for
% itself, and others acting on its behalf, a paid-up nonexclusive,
% irrevocable worldwide license in said article to reproduce, prepare
% derivative works, distribute copies to the public, and perform
% publicly and display publicly, by or on behalf of the Government. }}
% \normalsize 
% \end{flushright}
\appendix
\section{Computing gradient and Hessian information}
\label{sec:GradientHessian}

Here we derive the expressions for evaluating gradient and Hessian information for the deterministic inverse problem described in \eqref{eqn:laplace_approximation}. Define 
\begin{align}
 \mathcal{J}(\x) := {\frac{1}{2} \| \overline{\y}_i - \mathcal{G}(\x, t_i) \|_{\Gamma_i^{-1}}^2} + \displaystyle {\frac{\tau}{2}  \left \lVert \x \right \rVert^2_{\Sigma^{-1}}}.
\end{align}

{\bf Gradient:} By the chain rule we have
\begin{align}\label{eqn:gradient}
        \displaystyle \nabla_{\x}\mathcal{J} = \underbrace{\left(\frac{\partial \mathcal{G}(\x, t_i)}{\partial \x}\right)^{\top}\Gamma^{-1} \left(\mathcal{G}(\x, t_i) -  \overline{\y}_i  \right)}_{\textrm{Data misfit}} + \underbrace{\tau \Sigma^{-1}(\x)}_{\textrm{Regularization}}. 
\end{align}
Evaluating the ``Regularization'' term is straightforward. The ``Data misfit'' term requires evaluating the adjoint sensitivities. This can be accomplished by solving the following adjoint system backwards:
\begin{align}\label{eqn:adjoint}
        \displaystyle \frac{\mathrm{d}\x^*}{\mathrm{d}t} = \left(\frac{\partial \mathcal{G}(\x, t_i)}{\partial \x}\right)^{\top}\x^*\,, \quad t \in [t_i, 0]\,,\quad  \x^* (t_i)  = \Gamma^{-1} \left(\mathcal{G}(\x, t_i) -  \overline{\y}_i  \right)\,.
\end{align}
To solve the adjoint equation \eqref{eqn:adjoint} requires the Fr\'echet derivative along the forward trajectory $([0, t_i])$. Hence, one can checkpoint the forward trajectory and then propagate the adjoint trajectory backwards to obtain the adjoint sensitivities.
The gradient then can be evaluated as
\begin{align}
        \displaystyle \nabla_{\x}\mathcal{J} = \x^*(0) + \tau \Sigma^{-1}(\x)\,.   
\end{align}

{\bf Hessian:} To compute $\nabla^2_{\x,\x}\mathcal{J} \mathbf{v}$, we compute the directional derivative of $\nabla_{\x}\mathcal{J}$ in the direction $\mathbf{v}$. Applying the chain rule to equation \eqref{eqn:gradient}, we obtain
\begin{align}\label{eqn:hessian_vector}
       \displaystyle \nabla^2_{\x,\x}\mathcal{J} \mathbf{v} = \left(\frac{\partial \mathcal{G}(\x, t_i)}{\partial \x}\right)^{\top}\Gamma^{-1}\frac{\partial \mathcal{G}(\x, t_i)}{\partial \x} \mathbf{v} + \left(\frac{\partial^2 \mathcal{G}(\x, t_i)}{\partial \x^2}  \mathbf{v}\right)^{\top} \Gamma^{-1} \left(\mathcal{G}(\x, t_i) -  \overline{\y}_i  \right)  + \Sigma^{-1} \mathbf{v} \,.
\end{align}
In practice, the Hessian-vector product is approximated by the Gauss-Newton Hessian-vector product:
\begin{align}\label{eqn:GNhessian_vector}
        \displaystyle \nabla^2_{\x,\x}\mathcal{J}^{\rm GN} \mathbf{v} = \left(\frac{\partial \mathcal{G}(\x, t_i)}{\partial \x}\right)^{\top}\Gamma^{-1}\frac{\partial \mathcal{G}(\x, t_i)}{\partial \x} \mathbf{v} +   \Sigma^{-1} \mathbf{v} \,.
 \end{align}
 Evaluating the second term of the right-hand side in this equation is straightforward. In order to evaluate the first term of the right-hand side, both Tangent linear sensitivities and adjoint sensitivities need to be evaluated. That is, we first solve the following tangent linear system,
 \begin{align}\label{eqn:tlm}
        \displaystyle \frac{\mathrm{d}\delta \x}{\mathrm{d}t} = \frac{\partial \mathcal{G}(\x, t_i)}{\partial \x} \delta \x\,, \quad t\in[0,t_i]\,, \quad \delta \x(0) = \mathbf{v}\,,
 \end{align}
 and then solve the adjoint equation in \eqref{eqn:adjoint} with the initial forcing as $\x^*(t_i) = \Gamma^{-1}\delta \x(t_i)$. The system \eqref{eqn:tlm} should be solved along with the forward model. We now are ready to evaluate the Hessian-vector product:
 \begin{align}
        \displaystyle \nabla^2_{\x, \x}\mathcal{J}^{\rm GN}\mathbf{v} = \x^*(0) + \tau \Sigma^{-1}(\mathbf{v})\,.   
\end{align}
 \end{document}